\documentclass[sn-mathphys,Numbered]{sn-jnl}% Math and Physical Sciences Reference Style
\usepackage{graphicx}%
\usepackage{multirow}%
\usepackage{amsmath,amssymb,amsfonts}%
\usepackage{amsthm}%
\usepackage{mathrsfs}%
\usepackage[title]{appendix}%
\usepackage{xcolor}%
\usepackage{textcomp}%
\usepackage{manyfoot}%
\usepackage{booktabs}%
\usepackage{algorithm}%
\usepackage{algorithmicx}%
\usepackage{algpseudocode}%
\usepackage{listings}%
\usepackage{enumerate}
\usepackage{setspace}
\usepackage{lineno}

\theoremstyle{thmstyleone}%
\newtheorem{theorem}{Theorem}[section]%
\newtheorem{lemma}{Lemma}[section]%

\theoremstyle{thmstyletwo}%
\newtheorem{remark}{Remark}[section]%

\theoremstyle{thmstylethree}%
\newtheorem{definition}{Definition}[section]%

\theoremstyle{thmstylefour}%
\newtheorem{corollary}{Corollary}[section]%

\numberwithin{equation}{section}
\begin{document}

	\title[Li--Yau estimates for a Finslerian nonlinear parabolic equation]{Li-Yau estimates for a nonlinear parabolic equation on Finsler metric measure spaces}
	

\author*[1]{\fnm{{Bin}} \sur{Shen}}\email{shenbin@seu.edu.cn}

\author[1]{\fnm{Yuhan} \sur{Zhu}}\email{yuhanzhu@seu.edu.cn}
\equalcont{These authors contributed equally to this work.}



\affil*[1]{\orgdiv{School of Mathematics}, \orgname{Southeast University}, \orgaddress{\street{Dongnandaxue Road 2}, \city{Nanjing}, \postcode{211189}, \state{Jiangsu}, \country{China}}}

%	\author{\fnm{Bin Shen and Dingli Xia}}


\abstract{In this manuscript, we explore the positive solutions to the
	Finslerian nonlinear equation
	\textcolor{black}{\begin{eqnarray*}\frac{\partial u}{\partial t} = \Delta^{\nabla u} u + au\log u + bu,\end{eqnarray*}}which is related to Ricci solitons and serves as the Euler-Lagrange equation to the Finslerian log-energy functional. We then obtain the global gradient estimate of its positive solution on a compact Finsler metric measure space with the weighted Ricci curvature bounded below. Furthermore, applying a more general method based on a new comparison theorem  {established} by the first author, we derive the local gradient estimate on non-compact forward complete Finsler metric measure spaces with the mixed weighted Ricci curvature bounded below, under additional assumption of finite bounds of misalignment and some non-Riemannian curvatures. Lastly, we prove the Harnack inequalities {for} such solutions, as well as boundness or gap property of global solutions.}




\maketitle

\section{Introduction}
Log-Sobolev inequalities were first introduced and studied by L. Gross \cite{gross1975logarithmic} for Gaussian probability measure. More recently, S. Ohta \cite{ohta2017nonlinear} and S. Yin \cite{yin2021some} proved that in the Finsler measure space $(M,F,\mu)$,  
\begin{align}\label{logsobo}
	 {C}\int_M f^2 \log f^2 d \mu \leq  \int_M F^2(\nabla f)d\mu
\end{align}
for any $f\in C^{\infty}_0(M)$ such that $\int_M f^2 d\mu =\operatorname{Vol}_{d\mu} (M)$, where $C$ is a constant that  depends on different curvature conditions. In order to find the sharp constant in (\ref{logsobo}), as F. R. K. Chung and S.-T. Yau \cite{chung1996logarithmic} did in the Riemannian case, we can also define the \textit{Finslerian log-Sobolev constant} by
\begin{align}\label{cfls}
	C_{FLS} = \inf_{f \neq 0}\frac{\int_M F^2(\nabla f)d\mu}{\int_M f^2 \log f^2 d \mu},
\end{align}
{so that $C\leq C_{FLS}$.} Assume that $u$ is a minimizer of the functional in (\ref{cfls}), then by the variational calculation, $u$ satisfies the Euler-Lagrange equation that
\begin{align}
	\Delta^{\nabla u} u + C_{FLS} u\log u^2 + C_{FLS} u= 0,
\end{align}
which inspires us to consider this type of nonlinear parabolic equation
\begin{align}\label{maineq}
	\frac{\partial u}{\partial t} = \Delta^{\nabla u} u + au\log u + bu,
\end{align}
where $a$ and $b$ are two real constants, $\Delta^{\nabla u}$ is nonlinear Finslerian Laplacian. 

Another motivation is to understand the Ricci soliton defined by R. Hamilton \cite{hamilton1995formation}, which is a Riemannian manifold $(M,g)$ with a smooth
function $f\in C^\infty(M)$, such that for some constant $c\in \mathbb{R} $, it satisfies 
\begin{align}\label{riccisoli}
	\operatorname{Ric} = cg + \nabla^2 f ~~~~~~~~\text{on~} M,
\end{align}
where $\nabla^2$ and $\operatorname{Ric}$ are the Hessian tensor and Ricci tensor of the metric $g$, respectively.  
The Finslerian Ricci solitons have been investigated by {H. Zhu \cite{zhu2022class}}. 
%and some other results of rigidity and characterizations can be seen  {in \cite{Garofalo-Mondino2014,zhu2023rigidity,cheng2023characterizations}}.

Setting $u = e^f$ in (\ref{riccisoli}), one can deduce that
\begin{align}
	\Delta u + 2cu \log u = (A_0 - nc)u,
\end{align}
where $A_0$ is a constant (see \cite{ma2006gradient}). For this kind of elliptic equation, L. Ma \cite{ma2006gradient} established a Li-Yau estimate on a complete Riemannian manifold with Ricci curvature bounded below by a negative constant. Afterwards, Y. Yang \cite{yang2008gradient} generalized the estimate to the parabolic equation $\partial_t u = \Delta u + au\log u + bu$. Then G. Huang and B. Ma \cite{huang2010gradient} also considered the weighted Laplacian $\Delta_f = \Delta - {\nabla f \cdot \nabla}$, and replaced curvature condition with the $N$-Bakry-{\'Emery} Ricci tensor. {More generally, for metric measure space $(X,d,\mathfrak{m})$ with $RCD^*(K,N)$ condition, N. Garofalo and A. Mondino \cite{Garofalo-Mondino2014} have established Li-Yau estimate and Harnack inequality.} Apart from Li-Yau {estimates}, some other types of gradient estimates like Souplet-Zhang {estimates} have been obtained on this nonlinear equation by {X. Fu \cite{fu2022gradient} and N. Dung \cite{dung2019sharp}}.  

Although numerous studies have been carried out on gradient estimates for equations on Riemannian manifolds or Riemannian metric measure spaces, there are few published results on Finsler manifolds yet. S. Ohta and K.-T. Sturm have introduced the weighted Ricci curvature on a Finsler {metric measure space} as a generalization of the Bakry-{\'Emery}  Ricci curvature in Riemannian geometry (cf. \cite{Ohta2008HeatFO} and \cite{ohta2009finsler}) and proved a Finslerian version of Bochner–Weitzenböck formula (cf. \cite{ohta2014bochner}), which they applied to derive the Li-Yau estimate for heat equation on a compact Finsler  {metric measure space}. However, no previous study has investigated the  {parabolic equation with logarithmic nonlinear term like (\ref{maineq}) on any Finsler metric measure space}.  {W}e first generalized the results in \cite{ma2006gradient}, \cite{yang2008gradient} and \cite{huang2010gradient} to the nonlinear parabolic equation (\ref{maineq}) on the Finsler  {metric measure space}. The first main result in compact Finslerian case is the following.
\begin{theorem}\label{mainthm_cpt}
	Let $(M,F,\mu)$ be an $n$-dimensional compact Finsler  {metric measure space} . For some integer $N \in (n,\infty)$, suppose the weighted Ricci curvature $Ric^N \geq -K$ where $K$ is a nonnegative constant. If $u(x, t)$ is a positive global solution to equation (\ref{maineq}) on $M \times [0, +\infty)$ where $a \neq 0$, let $f(x, t) = \log u(x, t)$ and $\alpha = \max \{a/2, -a/4\}$. Then, for any $\beta > 1$, $x\in M$ and $t>0$, we have
	\begin{align}\label{ge_cpt}
		F^2(\nabla f)(x,t) + \beta af(x,t) + \beta b - \beta f_t(x,t) \leq \frac{N\beta^2}{2}\left\{\frac{1}{t} + \alpha + \frac{K}{2(\beta -1)}\right\}.
	\end{align}
\end{theorem}
In particular, as $a  { = } 0$ and $b=0$ , the equation (\ref{maineq}) simplifies to a common heat equation on a Finsler manifold. In this special case, the gradient estimate (\ref{ge_cpt}) would be the same as the result by Ohta and Sturm (cf. Theorem 4.4 in \cite{ohta2014bochner}). It is remarkable that S. Ohta \cite{ohta2009finsler} showed that the lower curvature bound ($\operatorname{Ric}^N \geqslant -K$) is equivalent to the
curvature-dimension condition $CD(-K, N)$ defined by  {Lott-Villani-Sturm \cite{lott2009ricci,sturm2006geometry1,sturm2006geometry2}.}

Further challenge, however, arose in the case of complete but non-compact Finsler manifolds. Although Q. Xia \cite{xia2023li} gave a gradient estimate on the non-compact case, that method is feasible only for Finslerian heat equation. Currently, the latest progress in comparison theorem made in \cite{shenac} filled the gap, so that we provided a systematic method to demonstrate the following gradient estimate.
\begin{theorem}\label{mainthm_cplt}
	Let $(M, F, \mu)$ be a forward complete $n$-dimensional Finsler metric measure space without boundary. Denote by $B_p(2R)$ the forward geodesic ball centering $p\in M$ with radius $2R$, and by $r(x):=d(p,x)$ the distance function from $p$.   {Suppose for some $N \in (n,\infty)$, the mixed weighted Ricci curvature $^m\operatorname{Ric}_{\nabla r}^N$ is bounded below by constant $-K(2R)$ and the  misalignment $\alpha \leqslant A(2R)$ in $B_p(2R)$, for $K(2R), A(2R)\geq 0$. Suppose the non-Riemannian tensors $U$, $\mathcal{T}$ (defined in (\ref{nonRiem_U}) and (\ref{nonRiem_T}) respectively) and $\operatorname{div}C(V):= C^{ij}_{~~k|i}(V)V^k\frac{\delta }{\delta x^j}$ of $M$
	satisfy the norm bound by 
	$$F(U) + F^*(\mathcal{T} ) + F^*(\operatorname{div}C) \leqslant K_0.$$}
	If $u(x, t)$ is a positive local solution to equation (\ref{maineq}) on $M \times [0, +\infty)$ where $a \neq 0$, by setting $f(x, t) = \log u(x, t)$ and $\alpha = \max \{a/2, -a/4\}$, then, for any $\beta > 1$, $0 <\delta< 1 $, $x\in B_p(R)$ and $t>0$, we have
	\begin{align}\label{ge_cplt}
		\begin{gathered}
			F^2(\nabla f)(x,t) + \beta af(x,t) + \beta b - \beta f_t(x,t) \leq \frac{N\beta^2}{2\delta}\left\{\frac{1}{t} + \alpha + \frac{AK}{2(\beta -1)}\right.\\
			\left.+ \frac{N\beta ^2 c_1^2}{16(1-\delta)(\beta - 1) R^2} + \frac{(2c_1^2 + c_2) A }{2R^2}+\frac{c_3(N,A,K_0)}{2R^2}(1+R+R\sqrt{K}) \right\},
		\end{gathered}
	\end{align}
	where $c_1$, $c_2$ and $c_3(N,A,K_0)$ are all positive constants.
\end{theorem} 

When dealing with non-compactness, the curvature conditions necessitate the bounds of misalignment and non-Riemannian quantities (further elaborated in Section 2 or \cite{shenac}). Our conditions are marginally stronger in contrast to those presented in \cite{xia2023li}. Nevertheless, such method is particularly useful in studying complete non-compact Finsler manifolds, and  {is expected to  apply to other nonlinear elliptic or parabolic equations on Finsler metric measure spaces, such as the parabolic Schr\"odinger equation \cite{shenac} and the elliptic Allen-Cahn equation \cite{shenAC2023}.} It is also worth noting that if the manifold degenerates to a Riemannian one, these conditions are just equivalent to the $(N-n)$-Bakry- {\'Emery} Ricci curvature bounded from below. Hence the gradient estimate (\ref{ge_cplt}) is a natural extension of the similar results in \cite{yang2008gradient} and \cite{huang2010gradient}.

Based on the aforementioned estimates, there are several noteworthy results. We will give a more general Harnack's inequalities than \cite{ohta2014bochner}.
\begin{corollary} Under the same conditions in Theorem \ref{mainthm_cpt}, we have
	\begin{align}\label{hi_cpt}
		[u(x_1, t_2)]^{1/{e^{a t_1}}} \leqslant[u(x_2, t_2)]^{{1}/{e^{a t_2}}}\exp (\Theta_\beta(x_1,x_2,t_1,t_2)),
	\end{align}
	where 
	\begin{align*}
		 {\Theta_\beta(x_1,x_2,t_1,t_2) =\int_{t_1}^{t_2} \frac{N\beta}{2e^{a\tau}}\left\{\frac{1}{\tau } + \alpha + \frac{K}{2(\beta -1)} + \frac{d(x_2,x_1)^2}{2N(t_2-t_1)^2}- \frac{2b}{\beta N}\right\} d\tau,}
	\end{align*}
	for any $\beta>1$, $0<t_1<t_2<T$ and $x_1$, $x_2\in M$.   
\end{corollary}

Taking $a \to 0$ and $b=0$ in (\ref{hi_cpt}), one will find that the Harnack inequality derived from the heat equation in \cite{ohta2014bochner} is also a special case of (\ref{hi_cpt}). Moreover, it can also be extended to non-compact Finsler  {metric measure space}.

\begin{corollary} Under the assumption in Theorem \ref{mainthm_cplt}, we have
	\begin{align}\label{hi_cplt}
		[u(x_1, t_2)]^{1/{e^{a t_1}}} \leqslant[u(x_2, t_2)]^{{1}/{e^{a t_2}}}\exp (\Theta_\beta(x_1,x_2,t_1,t_2)),
	\end{align}
	where 
	\begin{align*}
		\Theta_\beta(x_1,x_2,t_1,t_2) =\int_{t_1}^{t_2} \frac{N\beta}{2e^{a\tau}}\left\{\frac{1}{\tau } + \alpha + \frac{AK}{2(\beta -1)} + \frac{d(x_2,x_1)^2}{2N(t_2-t_1)^2}- \frac{2b}{\beta N}\right.\\
		\left.+ \frac{N\beta ^2 c_1^2}{16(1-\delta)(\beta - 1) R^2} + \frac{(2c_1^2 + c_2) A }{2R^2}+\frac{c_3(N,A,K_0)}{2R^2}(1+R+R\sqrt{K}) \right\} d\tau,
	\end{align*}
	for any $0 <\delta< 1$, $\beta>1$, $0<t_1<t_2<T$ and $x_1$, $x_2\in B_p(R)$ .
\end{corollary}

In particular, assume that the solution $u$ is independent of $t$, that is, $u$ solves the elliptic equation
\begin{align}\label{e_equa}
	\Delta^{\nabla u} u + au\log u + bu = 0 \text{~~on~~} M,
\end{align}
where $a \neq 0$ {.} We then have the following  {interesting property for global solution.}
\begin{theorem}\label{Liou}
	Let $(M, F, \mu)$ be a forward complete $n$-dimensional Finsler metric measure space without boundary, with misalignment bounded by $A$ {, and  these non-Riemannian tensors $U$, $\mathcal{T}$ and $\operatorname{div}C(V)$
	satisfy the norm bound by $$F(U) + F^*(\mathcal{T} ) + F^*(\operatorname{div}C) \leqslant K_0.$$} Let $u$ be a positive  {global} solution to (\ref{e_equa}) on $M$ and $f = \log u$. Suppose mixed weighted Ricci curvature  {$^m\operatorname{Ric}_{\nabla r}^N \geqslant -K$} for some $N$ and $K $ satisfying
	\begin{align}\label{Liou_con}
		n<N<\frac{2b}{\alpha}~~~~\text{and}~~~~  0\leqslant AK\leqslant\frac{4b}{N}+2\alpha-4\sqrt{\frac{2b\alpha}{N}},
	\end{align}
	where $\alpha = \max \{a/2, -a/4\}$. 

	\noindent
	 {\textbf{Case 1:} when $a<0$, then $u(x)\geqslant 1$ on $M$.}

	\noindent
	 {\textbf{Case 2:} when $a>0$, then $u(x)\leqslant 1$ on $M$.}
\end{theorem}
 {In case of $b=0$ and $K=0$, the property of global solution is shown in the theorem as follows.
\begin{theorem}\label{Lioub0k0}
	Let $(M, F, \mu)$ be a forward complete $n$-dimensional Finsler metric measure space without boundary. Suppose that $M$ has nonnegative mixed weighted Ricci curvature $^m\operatorname{Ric}_{\nabla r}^N \geqslant 0$ for some $N$, and misalignment is bounded by $A$. Also the non-Riemannian tensors $U$, $\mathcal{T}$ and $\operatorname{div}C(V)$
	satisfy $$F(U) + F^*(\mathcal{T} ) + F^*(\operatorname{div}C) \leqslant K_0.$$ Let $u$ be a positive global solution to
	\begin{align}\label{e_equa_b0}
		\Delta^{\nabla u} u + au\log u = 0 \text{~~on~~} M.
	\end{align} 
	\noindent
	\textbf{Case 1:} when $a<0$, then $u(x)\geqslant e^{-N/8}$ on $M$.
	\\
	\noindent
	\textbf{Case 2:} when $a>0$, then $u(x)\leqslant e^{N/4}$ on $M$.
\end{theorem}}

 {The previous research focused on the case of} $b=0$ and complete non-compact Riemannian manifold with nonnegative Ricci curvature. Y. Yang \cite{yang2008gradient} has proved that for $a > 0$, $u(x) \leqslant e^{n/2}$, and for $a < 0,  {u}(x) \geqslant e^{-n/4}$, indicating that $u(x)$ has a gap when $a<0$. Then X. Cao, et al. \cite{CAO20132312} further improved the gap to $u(x) \geqslant e^{-n/16}$ when $M$ is non-compact and $u(x) \geqslant 1$ when $M$ is closed.   {Theorem \ref{Lioub0k0} improves the outcome in  \cite{yang2008gradient} by providing a sharper boundness. Moreover, Theorem \ref{Liou} above shows that when $b>n\alpha/2 $ and $K\geqslant 0$, $u(x)$ also has similar boundness property when $a>0$ and gap property when $a<0$.}

This manuscript is structured as follows. In Section 2, we offer a brief  {introduction to basic} notions in Finsler geometry. Furthermore, we introduce the Finslerian log-energy functional, the Finslerian log-Sobolev constant and related variational equation. Then in Section 3, we demonstrate global gradient estimates of positive solutions with the $CD(-K, N)$ condition. In Section 4, we show local gradient estimates with the mixed weighted Ricci curvature condition. Finally, Section 5 is devoted to presenting some applications such as Harnack inequalities and  {analyzing the properties of the global solutions}.
\section{Preliminaries}
In this section, we briefly review some definitions and theorems in Finsler geometry (see \cite{bao2000introduction} and \cite{shen2001lectures}).
\subsection{Finsler manifold, connection, curvatures and misalignment}
A Finsler manifold is a pair $(M, F)$, where $M$ is an $n$-dimensional connected smooth manifold and $F: T M \rightarrow[0, \infty)$, called \textit{Finsler structure}, is a nonnegative function satisfying
\begin{itemize}
	\item $F \in {C}^{\infty}(T M \backslash \mathbf{0})$, where $\mathbf{0}$ is the zero section;
	\item $F(x, cv)=c F(x,v)$ for all $(x,v) \in T M$ and $c \geqslant 0$
	\item The $n \times n$ matrix
	\begin{align}
		\left(g_{i j}(x,v)\right)_{i, j=1}^n:=\left(\frac{1}{2} \frac{\partial^2\left(F^2\right)}{\partial v^i \partial v^j}(x,v)\right)_{i, j=1}^n
	\end{align}
	is positive-definite for all $(x,v) \in TM \backslash \mathbf{0}$.
\end{itemize}

In term of Finsler structure $F$, any  {non-vanishing vector field $V$ on $M$ } induces a Riemannian structure $g_V$ by
\begin{align}
	g_V(X,Y) = g_{ij}(V)X^iY^i,~~ X,Y\in TM
\end{align}
and a norm $F_V(\cdot)$ by
\begin{align}
	F^{ {2}}_V(X) = g_V(X,X).
\end{align}
In particular, the norm of reference vector is $F_V(V) = F(V)$.

There is an almost $g$-compatible and torsion-free connection $D^V$ called \textit{Chern connection} on the pull-back tangent bundle $\pi^*TM$. Namely, for any $X$, $Y$, $Z\in TM$,
\begin{align}
	\begin{gathered}
		 {Z}\left(g_V(X, Y)\right)-g_V\left(D^V_Z X, Y\right)-g_V\left(X, D^V_Z Y\right)=2 C_V\left(D^V_Z V, X, Y\right),\\
		D^V_X Y-D^V_Y X=[X, Y],
	\end{gathered}
\end{align}
where the \textit{Cartan tensor} $C_V$ is one of non-Riemannian tensors given by
\begin{align}
	C_V(X,Y,Z):=C_{ijk}(V)X^iY^jZ^k = \frac{1}{4}\frac{\partial^3 F^2(V)}{\partial v^i \partial v^j\partial v^k }X^iY^jZ^k.
\end{align}
Also letting $\Gamma^i_{jk}(V)$ denote the coefficients of Chern connection, one defines the \textit{spray coefficients} as 
\begin{align}
	G^i := \frac{1}{2}\Gamma^i_{jk}v^{ {j}}v^k,
\end{align}
and the \textit{spray} 
\begin{align}
	G =v^i\frac{\delta}{\delta x^i}= v^i\frac{\partial}{\partial x^i} - 2G^i\frac{\partial }{\partial v^i},
\end{align}
where $\frac{\delta}{\delta x^i} = \frac{\partial }{\partial x^i} - N^j_i\frac{\partial }{\partial v^j}$ and $N^j_i = \frac{\partial G^j}{\partial v^i}$.
Customarily, we denote the horizontal and vertical Chern derivatives by ``$|$" and ``; ", respectively. For example, for some tensor $T = T_idx^i$ on the pull-back bundle,
$$
T_{i|j}:= \frac{\delta T_i}{\delta x^j} - \Gamma^k_{ij}T_k,~~~~T_{i;j}:=\frac{\partial T_i}{\partial v^j} .
$$

The spray coefficients  {can also induce} the \textit{geodesic}, a curve $\gamma$ on $M$ satisfying 
\begin{align}
	\ddot{\gamma}(t) = 2G^i(\gamma(t),\dot{\gamma}(t)).
\end{align}
A \textit{forward geodesic ball} centered at $p$ with
radius $R$ can be represented by 
\begin{align}
	B^+_p(R):=\{q\in M: d(p,q)<R\}.
\end{align}
Here the \textit{forward distance} from $p$ to $q$ is given by
\begin{align}
	d(p,q) := \inf_{\gamma}\int_{0}^{1}F(\gamma(t),\dot{\gamma}(t))dt,
\end{align}
where the infimum is taken over all the $C^1$
curves $\gamma$ : $[0, 1] \to  M$ such that $\gamma(0) = p$
and $\gamma(1) = q$. Different from Riemannian manifold, the distance function is not necessarily symmetric between $p$ and $q$.
Adopting the exponential map, a Finsler manifold $(M, F)$ is said to be \textit{forward complete} if the exponential map is defined over $TM$. Thus, any two points in a forward complete manifold $ M$ can be connected by a minimal geodesic. Additionally, the forward closed ball $\overline{B^+_p(R)}$ is compact.  {In this case, we define the forward distance function from fixed point $p$ by $r(x)=d(p,x)$.}

The curvature form $\Omega$ induced by Chern connection can be divided into \textit{Chern Riemannian curvature} $R$ and \textit{Chern non-Riemannian connection} $P$ as
\begin{align}
	\Omega(X,Y)Z = R(X,Y)Z + P(X,D^V_YV,Z).
\end{align}
Then the Riemannian and non-Riemannian tensors are
\begin{align}
R^{i}_{jkl}&=\frac{\delta \Gamma^i_{jl}}{\delta x^k} - \frac{\delta \Gamma^i_{jk}}{\delta x^l} + \Gamma^i_{km}\Gamma^m_{jl} - \Gamma^i_{lm}\Gamma^m_{jk},\\
P^{i}_{jkl} &= \frac{\partial\Gamma^i_{jk} }{\partial y^l}.
\end{align}
Denoting $R_{ijkl}:= R^{s}_{ikl}g_{sj}$, the \textit{flag curvature} with pole $v$ is locally expressed by
\begin{align}
	K(u, v):=\frac{-R_{i j k l}(v) u^i v^j u^k v^l}{\left(g_{i k}(v) g_{j l}(v)-g_{i l}(v) g_{j k}(v)\right) u^i v^j u^k v^l} ,
\end{align}
for any two linearly independent vectors $u, v \in T_x M \backslash\mathbf{0}$, which span a tangent plane $\Pi=\operatorname{span}\{u, v\}$. Then the \textit{Finslerian Ricci curvature} is defined by
\begin{align}
	\operatorname{Ric}(v):=F^2(v) \sum_{i=1}^{n-1} K\left(e_i,v \right),
\end{align}
where $\{e_1,...e_{n-1}, v/F(v)\}$ is an orthonormal basis of $T_xM$ with respect to $g_v$. The \textit{T-curvature} is another non-Riemannian quantity, given by
\begin{align}
	T_y(v):=g_y\left(D_v V, y\right)-\hat{g}_x\left(\hat{D}_v V, y\right),
\end{align}
where $v \in T_x M, V$ is a vector field with $V(x)=v$, and $\hat{D}$ denotes the Levi-Civita connection of the induced Riemannian metric $\hat{g}=g_Y$. The T-curvature vanishes if and only if the Chern non-Riemannian curvature $P$ vanishes. To find the trace of T-curvature, let $\left\{e_i\right\}$ be an orthonormal basis with respect to metric $g_V$ at point $x$, where $V$ is a fixed reference vector field. Moreover, let $\left\{E_i\right\}$ be the local vector fields obtained by moving $\left\{e_i\right\}$ in parallel in a neighborhood of $x$ on $M$.  {D}efine the tensor $U$ as
\begin{align}
	U_v(W)=g(x, W)(U(v, W), W),
\end{align}
for any local vector field $W$, with
\begin{align}\label{nonRiem_U}
	U(v, W)=\sum_{i=1}^n\left(D_{e_i}^W E_i-\hat{D}_{e_i} E_i\right).
\end{align}

We introduce an important constant on a Finsler manifold defined in \cite{shenac}, called the \textit{misalignment}.
\begin{definition}[misalignment]
	For a Finsler manifold $(M, F)$, the misalignment of a Finsler metric at $x\in M$ is defined by
	\begin{align}
		\alpha(x):=\sup _{V, W, Y \in S_x M} \frac{g_V(Y, Y)}{g_W(Y, Y)},
	\end{align}
	and the global misalignment of the Finsler metric by
	\begin{align}
		\alpha:=\sup _{x \in M} \sup _{V, W, Y \in S_x M} \frac{g_V(Y, Y)}{g_W(Y, Y)}.
	\end{align}
	 {where $S_x M:=\{y\in T_xM:F(x,y) = 1\}$ is the unit sphere bundle of $M$ at $x$.}
\end{definition}
In particular, a Finsler manifold $(M, F)$ is Riemannian if and only if $\alpha(x) \equiv 1$. More characterizations of the misalignment can be found in \cite{shenac}.
\subsection{Finslerian gradient, Laplacian and mixed weighted Ricci curvature}
Fix $x\in M$, then the \textit{Legendre transformation} $\mathcal{L}^*$ maps $\alpha \in T_x^* M$ to the unique element $v \in T_x M$ such that $\alpha(v)=F^*(\alpha)^2$, where $F^*(\alpha):=F(v)$ is the dual norm of $F$. For a differentiable function $f: M \rightarrow \mathbb{R}$, the \textit{gradient vector} of $u$ at $x$ is defined by $\nabla u(x):=\mathcal{L}^*(D u(x)) \in T_x M$. If $D u(x)=0$, then clearly $\nabla u(x)=0$. If $D u(x) \neq 0$,  we can write
$$
\nabla u=\sum_{i, j=1}^n g^{i j}(\nabla u) \frac{\partial u}{\partial x^j} \frac{\partial}{\partial x^i}.
$$
On $M_u:= \{x\in M : D u(x) \neq\mathbf{0}\}$, we also define the \textit{Hessian} of $u$ by
$$
\operatorname{Hess} u(X,Y) := g_{\nabla u}(D_X^{\nabla u}\nabla u, Y) = g_{\nabla u}(D_Y^{\nabla u}\nabla u, X).
$$
The  {symmetry} of Hessian tensor can be shown in \cite{ohta2014bochner}.

A Finsler metric measure space $(M, F, \mu)$ is a Finsler manifold equipped with  a given measure $\mu$.  Assume $\mu$ to be absolutely continuous with respect to the Lebesgue measure on every chart with smooth positive density. Hence in local coordinates $\left\{x^i\right\}_{i=1}^n$, the volume form can be expressed as $d \mu= {e^{\Phi}} d x^1 \wedge \cdots \wedge d x^n$  {where $\Phi(x)$ is a smooth function on $M$}. For any $y \in T_x M \backslash\{0\}$, define
the \textit{distortion} by
$$
 {\tau(x, y):=\log {\sqrt{\operatorname{det} g_{i j}(x, v)}} - \Phi(x)}.
$$

For any point $x \in M$, let $\gamma=\gamma(t)$ be a forward geodesic from $x$ with the initial tangent vector $\dot{\gamma}(0)=v$. The \textit{$S$-curvature} of $(M, F, \mu)$ is
\begin{align}
	S(x, v):=\frac{d}{d t} \tau=\left.\frac{d}{d t}\left(\frac{1}{2} \log \operatorname{det}\left(g_{i j}\right)- {\Phi(x)}\right)(\gamma(t), \dot{\gamma}(t))\right|_{t=0}.
\end{align}
Modeling the definition of T-curvature,  {the first author} \cite{shenac} defined another non-Riemannian curvature
\begin{align}\label{nonRiem_T}
	\mathcal{T}(V, W):=\nabla^V \tau(V)-\nabla^W \tau(W),
\end{align}
for vector fields $V$, $W$ on $M$. Intuitively, $S$-curvature is the changing of distortion along the geodesic in direction $v$ and $\mathcal{T}$ is the difference of $\nabla \tau$ on the tangent sphere.

In local coordinates $\{x^i\}$, let $d\mu  = e^{\Phi}dx^1...dx^n$ express the volume form, then the \textit{divergence} of a smooth vector field $ V$ can be written as
\begin{align}
	\operatorname{div}_{d\mu}(V):= \sum_{i=1}^{n}\left(\frac{\partial V^i}{\partial x^i} + V^i\frac{\partial \Phi }{\partial x^i}\right). 
\end{align}
The \textit{Finslerian Laplacian} of a function $f$ on $M$  {is} given by
\begin{align}
	\Delta_{d\mu} f := \operatorname{div}_{d\mu}(\nabla f). 
\end{align}
On $M_f := \{x\in M : \nabla f \neq \mathbf{0}\}$, one notices that $\Delta_{d\mu}^{\nabla f}f = \Delta_{d\mu}f$ 
where $\Delta_{d\mu}^{\nabla f} $ is the weighted Laplacian defined by $\Delta_{d\mu}^{\nabla f} f:=\operatorname{div}_{d\mu}\left(\nabla^{\nabla f} f\right)$
where
\begin{align}
	\nabla^{V} f:=\left\{\begin{array}{ll}
		\sum_{i, j=1}^n g^{i j}(V) \frac{\partial f}{\partial x^j} \frac{\partial}{\partial x^i} & \text { on } M_f ,\\
		0 & \text { on } M \backslash M_f,
	\end{array}  
	\right.
\end{align}
for a non-vanishing vector field $V$. In particular, $\nabla^{\nabla f}f = \nabla f$ also holds on $M_f$. So it is valid to use $\nabla = \nabla^{\nabla f}$ and $\Delta = \Delta^{\nabla f}= \Delta^{\nabla f}_{d\mu}$ for short, when $f$ is exactly the function that $\nabla$ or $\Delta$ act on.
One also observes for later frequent use that, given functions $f$, $f_1$ and $f_2$ on $M_f$,
\begin{align}
	D f_2\left(\nabla^{\nabla f} f_1\right)=g_{\nabla f}\left(\nabla^{\nabla f} f_1, \nabla^{\nabla f} f_2\right)=D f_1\left(\nabla^{\nabla f} f_2\right).
\end{align}

Let $H^1(M) := W^{1,2}(M)$ be the Sobolev space under norm 
\begin{align}\label{nm}
	\Vert u \Vert_{H^1(M)}:=\Vert u\Vert_{L^2(M)} + \frac{1}{2}\Vert F(\nabla u)\Vert_{L^2(M)} \frac{1}{2} + \Vert \overleftarrow{F}(\overleftarrow{\nabla} u)\Vert_{L^2(M)},
\end{align}
where $\overleftarrow{F}(x,v) := F(x,-v)$ is \textit{the reverse Finsler metric} and $\overleftarrow{\nabla}$ is the gradient with respect to the reverse metric $\overleftarrow{F}$. Then $H^1(M)$ is in fact a Banach space with respect to $\Vert\cdot  \Vert_{H^1(M)}$, and $H_0^1(M)$  {denotes} the closure of $C^\infty_0(M)$ under the norm defined in (\ref{nm}).

Due to the lack of regularity, Finslerian Laplacian is usually viewed in distributional sense (or weak sense) that  
\begin{align}
	\int_M \varphi \Delta^V f d \mu=-\int_M D \varphi\left(\nabla^V f\right) d \mu,
\end{align}
 {for any function} $\varphi \in H_0^1(M) \cap L^{\infty}(M)$.

Employing $S$-curvature, we present the definition of the \textit{weighted Ricci curvature} introduced by S. Ohta \cite{ohta2014bochner} as follows.
\begin{definition}[\cite{ohta2014bochner}]
	Given a unit vector $v \in T_xM   $ and  {a positive number}
	$N \in [n,+\infty]$, the weighted Ricci curvature is defined by
	\begin{itemize}
		\item $\operatorname{Ric}^n(x,v):= \begin{cases}\operatorname{Ric}(x,v)+\dot{S}(x,v) & \text { if }{S}(x,v)=0; \\ -\infty & \text { if }{S}(x,v)\neq 0.\end{cases}$
		\item $\operatorname{Ric}^N(x,v):=\operatorname{Ric}(x,v)+\dot{S}(x,v) - \frac{S^2(x,v)}{N-n}$ when $n<N<\infty$.
		\item $\operatorname{Ric}^\infty (x,v):=\operatorname{Ric}(x,v)+\dot{S}(x,v)$.
	\end{itemize}
	Here $\dot{S}(x,v)$ is the derivative along the geodesic from $x$ in the direction of $v$.
\end{definition} 

B. Wu \cite{wu2015comparison} defined the \textit{weighted flag curvature} when $k\neq n$ and  {\cite{shenac}} completed this concept for any $k$.
\begin{definition}[\cite{wu2015comparison}\cite{shenac}]
	Given two linearly independent vectors $v$, $w\in T_xM$ and  {a positive number}
	$k \in [n,+\infty]$, the weighted Ricci curvature is defined by
	\begin{itemize}
		\item $\operatorname{K}^n(v,w):= \begin{cases}\operatorname{K}(v,w)+\frac{\dot{S}(x,v)}{(n-1)F^2(v)} & \text { if }{S}(x,v)=0; \\ -\infty & \text { if }{S}(x,v)\neq 0.\end{cases}$
		\item $\operatorname{K}^k(v,w):=\operatorname{K}(v,w)+\frac{\dot{S}(x,v)}{(n-1)F^2(v)} - \frac{S^2(x,v)}{(n-1)(k-n)F^2(v)}$ when $n<k<\infty$.
		\item $\operatorname{K}^\infty (x,v):=\operatorname{K}(v,w)+\frac{\dot{S}(x,v)}{(n-1)F^2(v)}$.
	\end{itemize}
\end{definition} 

Moreover,  {it has} also defined the \textit{mixed weighted Ricci curvature}.
\begin{definition}[\cite{shenac}]\label{shenac}
	Given two linearly independent vectors $v$, $w\in T_xM$ and 
	$k \in [n,+\infty]$, the mixed weighted Ricci curvature is defined by
	\begin{itemize}
		\item $^m\operatorname{Ric}^n_w(x,v):= \begin{cases}\operatorname{tr}_wR_v(v)+\dot{S}(x,v) & \text { if }{S}(x,v)=0 ;\\ -\infty & \text { if }{S}(x,v)\neq 0.\end{cases}$
		\item $^m\operatorname{Ric}^k_w(x,v):=\operatorname{tr}_wR_v(v)+\dot{S}(x,v) - \frac{S^2(x,v)}{k-n}$ when $n<N<\infty$.
		\item $^m\operatorname{Ric}^\infty_w (x,v):=\operatorname{tr}_wR_v(v)+\dot{S}(x,v)$.
	\end{itemize}
	where $\operatorname{tr}_wR_v(v) := g^{ij}(w)g_v(R_v(e_i,v), e_j)$ is the trace of flag curvature with respect to $g_w$.
\end{definition} 

 { We define the mixed weighted Ricci curvature $^mRic^k_{W}(V)$ to be the weighted Ricci curvature $Ric^k(V)$ when the reference vector $W=0$.
	%It looks like the mixed weighted curvature $^m\operatorname{Ric}^k_w(x,v)$ is defined only at the point where the vector field $w$ is non-zero. However, we can expand it when $w$ does not identically vanish. For a locally nontrivial vector field $w$, with the set $M_w:=\{x\in M \mid w(x)\neq0\}$ has zero measure. We can choose a orthonormal frame fields $\{e_1,\cdots,e_n\}$ with respect to $g_w=g_{ij}(x,w)dx^i\otimes dx^j$ wherever $w\neq 0$, so that the Ricci curvature part in Definition \ref{shenac} is given by
	%\begin{eqnarray}\label{mRicw}
	%	\operatorname{tr}_wR_v(v):=F^2(v) \sum_{i=1}^{n-1} K\left(e_i,v \right).
	%\end{eqnarray}
	%For the point where $w=0$, we take a limit on the orthogonal frame, and define the limitation of (\ref{mRicw}) as the value of $\operatorname{tr}_wR_v(v)$ in the mixed weighted Ricci curvature at $w=0$. Because the definition of $\operatorname{tr}_wR_v(v)$ as a tensor is independent of the choice of the frames, so this limit is well defined.
}
Note that the weighted Ricci curvature is a special case of the mixed weighted Ricci curvature since $^m\operatorname{Ric}^k_v(x,v) = \operatorname{Ric}^k(x,v)$.

 {We call the mixed weighted Ricci curvature bounded from below by a constant $K$, denoted by $^m\operatorname{Ric}^k_w\geq K$, if the inequality holds that $^m\operatorname{Ric}^k_w(x,v)\geq KF^2(v)$.}

Next we introduce two important tools, the Bochner inequality and the Laplacian comparison theorem.
\begin{theorem}[ {Bochner inequality} \cite{ohta2014bochner}]\label{othaBoch}
	Let $(M, F, \mu)$ be an $n$-dimensional Finsler metric measure space. Given $u \in H_{\mathrm{loc}}^2(M) \cap {C}^1(M)$ with $\Delta u \in H_{\mathrm{loc}}^1(M)$, for $N \in[n, +\infty]$, we have
	\begin{align}
		\Delta^{\nabla u}\left(\frac{F(\nabla u)^2}{2}\right)-D(\Delta u)(\nabla u) \geqslant \operatorname{Ric}^N(\nabla u)+\frac{(\Delta u)^2}{N}
	\end{align}
	in sense of distribution, that is, 
	\begin{align}
		-\int_M D \varphi\left(\nabla^{\nabla u}\left(\frac{F(\nabla u)^2}{2}\right)\right) d \mu \geqslant \int_M \varphi\left\{D({\Delta} u)(\nabla u)+\operatorname{Ric}^N(\nabla u)+\frac{({\Delta} u)^2}{N}\right\} d \mu
	\end{align}
	for all nonnegative functions $\varphi \in H_0^1(M) \cap L^{\infty}(M)$.
\end{theorem}
\begin{theorem}[Laplacian comparison theorem \cite{shenac}]\label{Shencom}
	Let $(M, F, \mu)$ be a forward complete $n$-dimensional Finsler metric measure space.  {Denote by $r(x)=d(p,x)$ the forward distance function start from $p\in M$, }and by $V$ a fixed vector field on $M$. Suppose the misalignment of $M$ is finite with upper bound $A$, and the mixed weighted Ricci curvature  {$^m\operatorname{Ric}^N_{\nabla r}$} of $M$ is bounded from below by $-K$ with $K>0$, for some $N>n$. Suppose the non-Riemannian curvatures $U$, $\mathcal{T}$ (defined in (\ref{nonRiem_U}) and (\ref{nonRiem_T}) respectively) and $\operatorname{div}C(V):= C^{ij}_{~~k|i}(V)V^k\frac{\delta }{\delta x^j}$
	satisfy the norm bound by $$F(U) + F^*(\mathcal{T} ) + F^*(\operatorname{div}C) \leqslant K_0.$$
	Then, by setting $l=K / C(N, A)$ with $C(N, A)=$ $N+(A-1) n-A$, wherever $r$ is $C^2$, the nonlinear Laplacian of $r$ with reference vector $V$ satisfies
	\begin{align}
		\Delta^V r \leq C(N, A) \sqrt{l}\coth (\sqrt{l}r)+\sqrt{A}K_0 .
	\end{align}
\end{theorem}
\subsection{Finslerian log-Sobolev constant and variational interpretation}
Let $(M, F, \mu )$ be a Finsler measure space. For any function $f\in H^1(M)$ with $\int_Mf^2d\mu=\operatorname{vol}_{d\mu }(M)$, define
the \textit{Finslerian log-energy functional} by
\begin{align}
	E(f):=\frac{\int_M F^2(\nabla f)d\mu}{\int_M f^2\log(f^2)d\mu },
\end{align}
and \textit{Finslerian log-Sobolev constant} $C_{FLS}$ is the infimum of $E$, we then have the following theorem:
\begin{theorem}
	If $u$ achieves log-Sobolev constant $C_{FLS}$ in the constraint of $\int_M u^2 d\mu = \operatorname{vol}_{d\mu}(M)$
	then the Euler-Lagrange equation is given by
	\begin{align}
		\Delta u + C_{FLS} u\log u^2 + C_{FLS} u = 0.
	\end{align}
\end{theorem}
\begin{proof}
	Choose any local function $v$ with compact support $M$ and $t\in \mathbb{R}$. By a direction variational calculation, we have
	\begin{align*}
		\begin{aligned}
			0=\frac{\partial }{\partial t}E(u+tv)_{t=0}=&\frac{2\int_M [D(v)(\nabla u) - C_{\nabla u}(\nabla u, \nabla u, \nabla v)]d\mu }{ \int_M u^2\log u^2 d\mu}\\
			& - \frac{\int_M F^2(\nabla u)d\mu \int_M v (2u\log u ^2+ 2u)d\mu }{\left(\int_M u^2\log u^2 d\mu\right)^2}.\\
		\end{aligned}
	\end{align*}
	 {By} the definition of $C_{FLS}$, it follows that
	\begin{align}
		\int_M D(v)(\nabla u) d\mu -C_{FLS}\int_M v (u\log u ^2+ u)d\mu=0.
	\end{align}
	The conclusion holds for the arbitrariness of $v$.
\end{proof}

 {Analogously} to \cite{Ohta2008HeatFO}, one can define the global and local solutions to equation (\ref{maineq}):
\begin{definition}[Global solution]
	For $T>0$, we say that a function $u$ on $ M\times[0, T] $ is a \textit{global solution} to the equation (\ref{maineq}) if $u \in L^2\left([0, T], H_0^1(M)\right) \cap H^1\left([0, T], H^{-1}(M)\right)$ and
	\begin{align}
		\int_M \varphi \frac{\partial u}{\partial t} d\mu=\int_M[ -D \varphi\left(\nabla u\right) + \varphi(au \log u + bu)]d\mu,
	\end{align}
	holds for all $t \in[0, T]$ and  {$\varphi \in C_0^{\infty}(M)$.}
\end{definition}
\begin{definition}[Local solution]
	Given an open interval $I \subset \mathbb{R}^+$ and an open set $\Omega \subset M$, we say that a function $u$  is a \textit{local solution} to the equation (\ref{maineq}) on $\Omega \times I$ if $u \in L_{\mathrm{loc}}^2(\Omega \times I), F(\nabla u) \in L_{\mathrm{loc}}^2(\Omega \times I)$ and
	\begin{align}
		\int_I \int_{\Omega} u \frac{\partial \varphi}{\partial t} d \mu d t=\int_I \int_{\Omega} [ D \varphi\left(\nabla u\right) - \varphi(au \log u + bu)] \mu d t,
	\end{align}
	holds for all  {$\varphi \in C_0^{\infty}(\Omega \times I)$}.
\end{definition}
\begin{remark}
	The elliptic and parabolic regularities in \cite{Ohta2008HeatFO} and \cite{shenac} guarantee that the global solution $u(x,t)$ of (\ref{maineq}) enjoys the $H^2_{\text{loc}}$-regularity in $x$ and $C^{1,\beta}$-regularity ($0<\beta<1$) in term of both $t$ and $x$, meanwhile $\partial_t u \in H^1_{\text{loc}}(M)\cap C(M)$.
	Moreover, for the local solution $u(x,t)$ on $\Omega\times  I$ then $u\in H^2(\Omega)\cap C^{1,\beta}(\Omega \times I)$, $\Delta u \in H^1(M)\cap C(\Omega)$ and $\partial_t u \in H^1_{\text{loc}}(\Omega)\cap C(\Omega)$.
\end{remark}

\section{Gradient estimate on compact Finsler  {metric measure space}}
In this section we will always assume that $(M, F, \mu)$ is a compact $n$-dimensional Finsler metric measure space, whose weighted Ricci curvature $\operatorname{Ric}^N$ has lower bound of $-K$ ($K\geqslant0$) for some $N \in (n,\infty)$.
Let  {$w$} be a positive solution to (\ref{maineq}) on $M \times [0,+\infty)$.  {Let $w = e^{-b/a}u $, then $\nabla w = e^{-b/a}\nabla u$  and it follows that} 
\begin{align}\label{maineq_b0}
	  {\frac{\partial u}{\partial t} = \Delta^{\nabla u} u + au\log u},
\end{align}i.e.
\begin{align}
	\int_{0}^{T}\int_{M}\left[D\varphi(\nabla u)  - u\frac{\partial \varphi}{\partial t} \right] d\mu dt= \int_{0}^{T}\int_{M}a\varphi u\log u d\mu dt,
\end{align}
for any test function  {$\varphi\in C_0^\infty(M\times [0,T])$}.

Before we begin the proof of Theorem \ref{mainthm_cpt}, we shall introduce some lemmas: 
\begin{lemma}\label{lem_u2f}
	For any fixed $T>0$, Let $u$ be a positive solution to (\ref{maineq_b0}) on $M\times [0,T]$ and set  {$f:= \log u$, and $f_t = \partial_t f$}, then
	\begin{enumerate}
		\item $ {-\Delta} f=  af + F^2(\nabla f) - f_t,$
		\item $ {-\Delta^{\nabla u}(tf_t) = 2D(tf_t)\left(\nabla f\right) + (f_t-af) - \partial_t\left(t(f_t-af)\right),}$
	\end{enumerate}
	in the distributional sense.
\end{lemma}
\begin{proof}
	Since $\nabla u = u \nabla f$ and $u>0$ we have $g_{\nabla u} = g_{\nabla f}$ and hence $\Delta^{\nabla u}f = \Delta^{\nabla f}f = \Delta f$ on $M_u = M_f$.   
	Then
	\begin{align*}
		\Delta f &= 	\operatorname{tr}_{\nabla f} (\operatorname{Hess} f) - S(\nabla f)\\
		&= \frac{\operatorname{tr}_{\nabla u}(\operatorname{Hess} f)}{u} - \frac{F^2(\nabla u)}{u^2} - S\left(\frac{\nabla u}{u}\right)\\
		&= \frac{\Delta u}{u} - F^2(\nabla f),
	\end{align*}
	in weak sense on $M_u$. Hence,
	\begin{align*}
		- \int_{0}^{T}\int_{M_u}\left[D\varphi(\nabla f) - \frac{\partial \varphi}{\partial t}f\right]d\mu dt &= \int_{0}^{T}\int_{M_u}\varphi\left[\Delta f - \frac{\partial f}{\partial t}\right]\\
		&= \int_{0}^{T}\int_{M_u}\varphi\left[\frac{1}{u}\left(\Delta u - \frac{\partial u}{\partial t}\right) - F^2(\nabla f)\right]d\mu dt\\
		&=\int_{0}^{T}\int_{M_u}\varphi\left( -af - F^2(\nabla f) \right) d\mu dt.
	\end{align*}
	Since $\Delta u= 0 $ a.e. on $M \backslash M_u$, 
	\begin{align}\label{eqf_w}
		\int_{0}^{T}\int_M\left[D\varphi(\nabla f) - \frac{\partial \varphi}{\partial t}f\right]d\mu dt =  \int_{0}^{T}\int_{M}\varphi\left( af + F^2(\nabla f) \right) d\mu dt.
	\end{align}
	 {for any test function $\varphi\in C_0^\infty(M\times [0,T])$}. It completes the first part.

	 {Similarly, on the other hand, 
	\begin{align}\label{eqf2_w}
			\int_{0}^{T}\int_{M}D\varphi&(\nabla^{\nabla u} tf_t)d\mu dt =\int_{0}^{T}\int_{M_u}D\varphi(\nabla^{\nabla u} tf_t)d\mu dt \notag\\
%			&= \int_{0}^{T}\int_{M_u}Dt\varphi(\nabla^{\nabla u} f_t)d\mu dt \notag\\ 
			&= -\int_{0}^{T}\int_{M_u}D\left(\partial_t\left(t\varphi\right)\right)(\nabla f)d\mu dt  \notag\\
			&= \int_{0}^{T}\int_{M_u}\partial_t\left(t\varphi\right)\left[f_t - af - F^2(\nabla f)\right]d\mu dt \notag\\
			&= \int_{0}^{T}\int_{M_u}\partial_t\left(t\varphi\right)\left(f_t - af \right) +t\varphi \partial_t\left( F^2(\nabla f)\right)d\mu dt \notag\\
			&= \int_{0}^{T}\int_{M_u}\varphi_tt\left(f_t - af \right) +\varphi\left(f_t - af \right) + 2\varphi D(tf_t)\left(\nabla f\right)d\mu dt \notag\\
			&= \int_{0}^{T}\int_{M}\varphi_tt\left(f_t - af \right) +\varphi\left(f_t - af \right) + 2\varphi D(tf_t)\left(\nabla f\right)d\mu dt,
	\end{align}
%	Similarly we obtain
%	\begin{align}
%	\int_{0}^{T}\int_{M}D\varphi(\nabla^{\nabla u} tf_t)d\mu dt = \int_{0}^{T}\int_{M}\varphi_tt\left(f_t - af \right) +\varphi\left(f_t - af \right) + 2\varphi D(tf_t)\left(\nabla f\right)d\mu dt
%	\end{align}
	which is exactly the second identity.}  
\end{proof}
\begin{lemma}\label{lem_f2H}
	Let $(M,F,\mu)$ be a Finsler metric measure space with weighted Ricci curvature $\operatorname{Ric}^N \geqslant -K$, where the constant $K\geqslant 0$. Let $f$ satisfy (\ref{eqf_w}) and  {set} $f_t = \partial_t f$. Then for any constant $\beta > 1$,
	\begin{align}\label{def_H}
		H:= t\left\{F^2(\nabla f) + \beta (af -f_t)\right\},
	\end{align}
	 {satisfies} that 
	\begin{align}\label{H_lowb}
		\int_{0}^{T}\int_M -D\varphi(\nabla^{\nabla u}H) + 2\varphi g_{\nabla u}\left(\nabla f , \nabla^{\nabla u} H\right) +\varphi_t H d\mu dt \geqslant \int_{0}^{T}\int_M \varphi \mathcal{J} d\mu dt,
	\end{align}
	for any  {nonnegative} test function  {$\varphi \in C_0^\infty( M\times[0,T]) $}, where
	\begin{align}\label{H_lowb_J}
		\mathcal{J}(x,t) := t\left\{ {\frac{2}{N}\left(-\frac{H}{\beta t}-\left(1-\frac{1}{\beta}\right) F^2(\nabla f)\right)^2}+((\beta-1) a-2 K) F^2(\nabla f)\right\}-a H-\frac{H}{t}.
	\end{align}
\end{lemma}
\begin{proof}
	By the same method in Lemma \ref{lem_u2f},  {it is valid to consider the integral over $M_u$. During the proof of Lemma \ref{lem_u2f}, we see
	$$\int_{0}^{T}\int_{M_u} D\varphi(\nabla f) d\mu dt =\int_{0}^{T}\int_{M_u} \varphi\left(af - f_t + F^2(\nabla f)\right) d\mu dt,$$}
	for any  {nonnegative} test function  {$\varphi \in C_0^\infty(M\times[0,T]) $}. From (\ref{def_H}), one obtains that {
		\begin{align}\label{H1}
			\int_{0}^{T}\int_{M_u} D\varphi(\nabla f) d\mu dt =& \int_{0}^{T}\int_{M_u}\varphi\left[ \frac{H}{\beta t} + (1-\frac{1}{\beta})F^2(\nabla f) \right]d\mu dt.
		\end{align}}
	
	%\begin{align}
	%	\Delta f =& -\frac{H}{\beta t} - (1-\frac{1}{\beta})F^2(\nabla f) \\
	%	=& -\frac{H}{t} + (\beta - 1)(af-f_t).
	%\end{align}

	%Note that $[\Delta^{\nabla u}, \partial_t]u = 0$ holds pointwise on $M_u$ (see \cite{xia2023li}), and  {differentiate} (3.3) and (3.5) with respect to $t$. It then follows that
	%\begin{align}
	%	\begin{aligned}
	%		\Delta^{\nabla u} f_t &=  \partial_t (\Delta^{\nabla u} f)  \\
	%		&=  - af_t - 2 g_{\nabla u}\left(\nabla f, \nabla^{\nabla u} f_t\right) + f_{tt},  
	%	\end{aligned}
	%\end{align}
	%	\begin{align}
	%	\begin{aligned}
	%   \int_{M} D\varphi(\nabla^{\nabla u} f_t) d\mu &= \int_{M} \partial_t D\varphi(\nabla^{\nabla u} f) d\mu \\
	%    &= \int_{M}\varphi\left[ af_t + 2 g_{\nabla u}\left(\nabla f, \nabla^{\nabla u} f_t\right) - f_{tt}   \right]d\mu
	%    \end{aligned}
	%	\end{align}
	%and
	%\begin{align}
	%	\begin{aligned}
	%		H_t & =F^2(\nabla f)+\beta\left(a f-f_t\right)+t\left\{2 g_{\nabla n}\left(\nabla f, \nabla^{\nabla u} f_t\right)+a \beta f_t-\beta f_{t t}\right\} \\
	%		& =\frac{H}{t}+t\left\{2 g_{\nabla u}\left(\nabla f, \nabla^{\nabla u} f_t\right)+a \beta f_t-\beta f_{tt}\right\}.
	%	\end{aligned}
	%\end{align}
	
	Now in order to estimate $\Delta^{\nabla u} H$ in weak sense, using (\ref{eqf2_w}) and (\ref{H1}), we have
	%\begin{align}
	%	\begin{aligned}
	%		\Delta^{\nabla u} H&=t\left\{\Delta^{\nabla u} F^2(\nabla f)+a\beta \Delta^{\nabla u} f-\beta \Delta^{\nabla u} f_t\right\} \\
	%		&\geqslant t\left\{2 \operatorname{Ric}^N(\nabla f)+2 g_{\nabla u}\left(\nabla f, \nabla^{\nabla u} \Delta f\right)+\frac{2(\Delta f)^2}{N}+a \beta \Delta^{\nabla u} f -\beta \Delta^{\nabla u} f_t\right\}.
	%	\end{aligned}
	%\end{align}
	 {
	\begin{align}\label{e1}
		\begin{aligned}
			&\int_{0}^{T}\int_{M_u} D\varphi(\nabla^{\nabla u} H) d\mu dt \\=& \int_{0}^{T}\int_{M_u}\left\{t D\varphi(\nabla^{\nabla u}F^2(\nabla f)) + a\beta tD\varphi(\nabla f) - \beta D\varphi(\nabla^{\nabla u} tf_t) \right\}d\mu dt\\
			=& \int_{0}^{T}\int_{M_u}\left\{t D\varphi(\nabla^{\nabla u}F^2(\nabla f)) + \left[aH + a(\beta - 1)tF^2(\nabla f)\right]\varphi - \beta \varphi_t\left( f_t -af \right) \right.\\
			&\quad\quad\quad\quad\,\,\left.-\beta \varphi\left( f_t -af \right) - 2\beta \varphi D(tf_t)(\nabla f)\right\}d\mu dt.
		\end{aligned}
	\end{align}
	By substituting (\ref{def_H}) to the right-hand side of (\ref{e1}), it infers that 
	\begin{align}\label{e2}
		\begin{aligned}
			&\int_{0}^{T}\int_{M_u} D\varphi(\nabla^{\nabla u} H) - H\varphi_t d\mu dt \\
			=& \int_{0}^{T}\int_{M_u}\bigg\{t D\varphi(\nabla^{\nabla u}F^2(\nabla f)) + \left[aH +\frac{H}{t} +  a(\beta - 1)tF^2(\nabla f)\right]\varphi\\
			& \quad\quad\quad\quad- \left(t\varphi_t + \varphi\right)F^2(\nabla f) - 2\beta \varphi D(tf_t)(\nabla f) \bigg\}d\mu dt\\
			=& \int_{0}^{T}\int_{M_u}\bigg\{t D\varphi(\nabla^{\nabla u}F^2(\nabla f)) + \left[aH +\frac{H}{t} + a(\beta - 1)tF^2(\nabla f)\right]\varphi\\
			& \quad\quad\quad\quad - 2\left(\beta -1\right) \varphi D(tf_t)(\nabla f) \bigg\}d\mu dt.
		\end{aligned}
	\end{align}}

	%\begin{align}
	%	\begin{aligned}
	%		\Delta^{\nabla u} H&=t\left\{\Delta^{\nabla u} F^2(\nabla f)+a\beta \Delta^{\nabla u} f-\beta \Delta^{\nabla u} f_t\right\} \\
	%		&\geqslant t\left\{2 \operatorname{Ric}^N(\nabla f)+2 g_{\nabla u}\left(\nabla f, \nabla^{\nabla u} \Delta f\right)+\frac{2(\Delta f)^2}{N}+a \beta \Delta^{\nabla u} f -\beta \Delta^{\nabla u} f_t\right\}.
	%	\end{aligned}
	%\end{align}

	We then need to utilize the  {Bochner inequality} in Theorem \ref{othaBoch}. Using  the assumption that $\operatorname{Ric}^N(\nabla f) \geqslant - K F^2(\nabla f)$, we have  {
	\begin{align}\label{e3}
		\begin{aligned}
			&\int_{0}^{T}\int_{M_u} t D\varphi(\nabla^{\nabla u}F^2(\nabla f)) d\mu dt \\
			\leqslant & \int_{0}^{T}\int_{M_u}t\left\{2 K F^2(\nabla f)  - \frac{2(\Delta f)^2}{N} - 2D(\Delta f)(\nabla f)\right\}\varphi d\mu dt\\
			=  & \int_{0}^{T}\int_{M_u}t\left\{2 K F^2(\nabla f)  - \frac{2(\Delta f)^2}{N} - 2D\left(f_t - af - F^2(\nabla f)\right)(\nabla f)\right\}\varphi d\mu dt\\
			=  & \int_{0}^{T}\int_{M_u}t\left\{2 K F^2(\nabla f)  - \frac{2(\Delta f)^2}{N} - 2D\left(-\frac{H}{t} + \left(\beta -1\right)\left(af-f_t\right)\right)(\nabla f)\right\}\varphi d\mu dt\\
			= & \int_{0}^{T}\int_{M_u}\bigg\{2 tK F^2(\nabla f)  - \frac{2t(\Delta f)^2}{N} + 2D\left(H\right)(\nabla f)-2a(\beta -1)tF^2(\nabla f)\\
			& \quad\quad\quad\quad + 2\left(\beta -1\right) \varphi D(tf_t)(\nabla f) \bigg\}d\varphi d\mu dt
		\end{aligned}
	\end{align}}
	%\begin{align}
	%	\begin{aligned}
	%		\Delta^{\nabla u} H&\geqslant t \left\{ -2 K F^2(\nabla f) + 2 g_{\nabla u}\left(\nabla f, \nabla^{\nabla u} \left(-\frac{H}{t} + (\beta - 1)(af-f_t) \right)\right)  +\frac{(\Delta f)^2}{N} \right.\\
	%		&~~~~+ a\beta\left.\left(-\frac{H}{t} + (\beta - 1)(af-f_t) \right) - \beta \left(f_{tt}- af_t - 2 g_{\nabla u}(\nabla f, \nabla^{\nabla u} f_t)\right) \right\} \\
	%		&= t \left\{ -\frac{2}{t} g_{\nabla u}(\nabla f, \nabla^{\nabla u}H) + \frac{(\Delta f)^2}{N} - \frac{a\beta}{t}H + \left(2(\beta-1)a - 2K \right)F^2(\nabla f) \right.\\
	%		&~~~~+  \left.a^2\beta\left(\beta-1\right)f - a\beta(\beta-2)f_t - \beta f_{tt} + 2g_{\nabla u}\left(\nabla f, \nabla^{\nabla u}f_t\right)\right\}.
	%	\end{aligned}
	%\end{align}
	
	Combining (\ref{e2}), (\ref{e3}) and (\ref{H1}) it implies that 
	\begin{align}
		\begin{aligned}
			\Delta^{\nabla u}H  + 2 g_{\nabla u}\left(\nabla f , \nabla^{\nabla u} H\right) - H_t\geqslant &\frac{2t}{N}\left(-\frac{H}{\beta t}-\left(1-\frac{1}{\beta}\right) F^2(\nabla f)\right)^2\\
			&+t((\beta-1) a-2 K) F^2(\nabla f)-a H-\frac{H}{t},
		\end{aligned}
	\end{align}
	in weak sense.
\end{proof}

Then we start to prove Theorem \ref{mainthm_cpt}.
\begin{proof}[Proof of Theorem \ref{mainthm_cpt}]
	Let $(z,s)$ be the maximal point of $H$ on $M \times[0,T]$ and $\Lambda := H(z,s)$ be the maximum of $H$. Without loss of generality, we can assume $\Lambda > 0$ and $s > 0$, otherwise, the assertion of the theorem obviously holds. Now we shall follow the proof in \cite{ohta2014bochner} by  {Ohta and Sturm}  to claim that $\mathcal{J}(z,s) \leqslant 0$. By the assumption of contrary, $\mathcal{J}(z,s) > 0$, it follows that $\mathcal{J} > 0$ on a neighborhood of $(z,s)$, say, a small parabolic cylinder $B_z(\delta) \times [s-\delta,s]$, where $B_z(\delta)$ is the forward geodesic ball with radius $\delta$. Thus, $H$ must be a strict subsolution to the linear parabolic operator  
	$$
	\operatorname{div}_\mu(\nabla^{\nabla u}H) + g_{\nabla u}\left(\nabla f , \nabla^{\nabla u} H\right) - \partial_t H .
	$$
	This implies that $H(z,s)$ would be strictly smaller than the supremum of $H$ on the boundary of the parabolic cylinder $B_z(\delta) \times [s-\delta,s]$, which yields contradiction to the maximality of $H(z,s)$. Therefore, at point $(z,s)$, from (\ref{H_lowb_J}) and (\ref{H1}), we obtain that 
	\begin{align}\label{after0}
		0 \geqslant s\left\{\frac{2}{N}\left(-\frac{\Lambda}{\beta s}-\left(1-\frac{1}{\beta}\right) F^2(\nabla f)\right)^2+((\beta-1) a-2 K) F^2(\nabla f)\right\}-a \Lambda-\frac{\Lambda}{s}.
	\end{align}
	Denoting $v:=F^2(\nabla f)/{H}$ and multiplying both sides of (\ref{after0}) by $s\Lambda$, we have  
	\begin{align}\label{H2v}
		a s+1 \geqslant \frac{2}{N \beta^2} \Lambda\left(1+(\beta-1)v s\right)^2+((\beta-1) a-2 K) v s^2.
	\end{align}
	
	\textit{Case 1: }Firstly we consider the case of $a<0$ which means $a s+1 \leqslant 1$, then (\ref{H2v}) becomes
	\begin{align}\label{lambd1}
		\Lambda \leqslant \frac{N \beta^2}{2} \cdot \frac{1+(-(\beta-1) a +2 K) s \cdot v s}{\left(1+(\beta-1) v s\right)^2}.
	\end{align}
	Notice that, for the constants $r,p \geqslant 0$ and $q > 0$, the inequality 
	\begin{align}\label{esti1}
		r+p x\leqslant \left(r+\frac{p}{4 q} \right)(1+q x)^2,
	\end{align}
	holds for any $x \geqslant 0$. Hence,
	\begin{align}
		\Lambda \leqslant \frac{N\beta^2}{2}\left\{1 - \frac{a}{4}s + \frac{K}{2(\beta - 1)}s\right\}.
	\end{align}
	
	\textit{Case 2:} On the other hand, when $a>0$, then $(\beta - 1)avs^2 \geqslant 0$, and similarly (\ref{H2v}) yields 
	\begin{align}\label{lambd2}
		\Lambda \leqslant \frac{N \beta^2}{2} \cdot \frac{1+a s+2 K v s^2}{(1+(\beta-1) v s)^2} \leqslant \frac{N \beta^2}{2}\left\{1+\frac{a}{2} s+\frac{K}{2(\beta-1)} s\right\}.
	\end{align}
	
	Finally, set $\alpha :=  \max \{a/2, -a/4\}$ and combine (\ref{lambd1}) with (\ref{lambd2}). Since 
	\begin{align}
		\begin{aligned}
			F^2(\nabla f)(x,T) + \beta af(x,T) - \beta f_t(x,T) &= \frac{H(x,T)}{T} \leqslant \frac{\Lambda}{T} \\
			&\leqslant \frac{N \beta^2}{2}\left\{\frac{1}{T}+\alpha \frac{s}{T}+\frac{K}{2(\beta-1)} \frac{s}{T}\right\}\\
			&\leqslant \frac{N \beta^2}{2}\left\{\frac{1}{T}+\alpha+\frac{K}{2(\beta-1)} \right\}.
		\end{aligned}
	\end{align}
	This completes the proof for any $T>0$.
\end{proof}
\section{Gradient estimate on non-compact Finsler  {metric measure space}}
In this section,  {we} will continue to discuss the case of complete but non-compact Finsler manifolds. Let us review the conditions. In forward complete $n$-dimensional Finsler metric measure space $(M, F, \mu)$, let $B_p(2R)$ be the forward geodesic ball centering $p \in M$ with radius $2R$, and $r(x):=d(p,x)$ be the distance function from $p$. Suppose for some $N \in (n,\infty)$, the mixed weighted Ricci curvature $^m\operatorname{Ric}_{\nabla r}^N \geqslant -K(2R)$, where $K(2R)\geqslant0$. Additionally, we need to suppose that $M$ has finite upper bounded misalignment and non-Riemannian tensors by $A(2R)$ and $K_0$ respectively. It is sufficient to consider $u(x, t)$ to be a positive solution to the equation (\ref{maineq_b0}) (i.e. $b=0$) on $B_p(2R) \times [0, T]$ where $a \neq 0$. Same as the method in section 3, we set $f(x, t) := \log u(x, t)$ and H:= $t\left\{F^2(\nabla f) + \beta (af -f_t)\right\}$ for any $\beta > 1$. The finite bounded misalignment also infers the weighted Ricci curvature $\operatorname{Ric}^N \geqslant -AK$.
Then Lemma \ref{lem_u2f} and Lemma \ref{lem_f2H} still hold on $B_p(2R) \times [0, T]$. Now we start the proof with choosing a cut-off function.
\begin{proof}[Proof of Theorem 1.2]
	We define a smooth function $\phi$ on $[0, \infty)$ such that $\phi(d) \equiv 1$ for $d \leqslant 1$, $\phi(d) \equiv 0$ for $d \geqslant 2$ and $0\leqslant \phi \leqslant 1$ otherwise. Assume that 
	\begin{align}\label{con_cut}
		-c_1 \leqslant \frac{\phi'(d)}{\phi^\frac{1}{2}(d)} \leqslant 0 \text{~~and~~}\phi''\geqslant-c_2
	\end{align}
	for two positive constants $c_1$ and $c_2$. Then we define the cut-off function on $M$ by
	\begin{align}
		\Psi(x) = \phi\left(\frac{r(x)}{R}\right) .
	\end{align}
	 {Clearly we notice that $\Psi(x)$ is supported in $B_p(2 R)$. However, in general $\Psi(x)$ is not smooth since $r(x)$ is only Lipschitz on the cut locus of $p$. Without loss of generality, we can assume the smoothness of $\Psi(x)$ by choosing the neighborhood of $p$ such that the intersection with the cut locus of $p$ is empty.} Then for any reference vector field $V$,
	\begin{align}\label{cut_F}
		\begin{aligned}
			F^2_{V}\left(\nabla^{V} \Psi\right) = \frac{1}{R^2}\left( \phi'(\frac{r}{R})\right)^2 F_{V}^2\left(\nabla^{V}r\right) \leqslant  \frac{c_1^2A}{R^2} \Psi,
		\end{aligned}
	\end{align}
	where the inequality is from (\ref{con_cut}) and the bound of misalignment:
	\begin{align}
		F_{V}^2(\nabla^{V} r) \leqslant A F^2_{\nabla^{V}r}\left(\nabla^{V}r\right) = A.
	\end{align}
	
	Applying the Laplacian comparison theorem in Theorem \ref{Shencom} and taking $V = \nabla u$ on $M_u$, we have 
	\begin{align}\label{cut_Delta}
		\begin{aligned}
			\Delta^{\nabla u} \Psi &=\phi^{\prime \prime}\left(\frac{r}{R}\right) \frac{F_{\nabla u}^2\left(\nabla^{\nabla u} r\right)}{R^2}+\phi^{\prime}\left(\frac{r}{R}\right) \frac{\Delta^{\nabla u} r}{R} \\
			&\geqslant-\frac{c_2 A}{R^2}-\frac{c_1}{R}\left(C(N, A) \sqrt{l} \operatorname{coth}(R \sqrt{l})+c_0\right)\\
			&\geqslant-\frac{c_2 A}{R^2}-\frac{c_1}{R^2}\left(C(N, A)\left(1+R\sqrt{l}\right)+c_0 R\right) .
		\end{aligned}
	\end{align}
	The last inequality comes from the inequality
	$
	\coth (x) \leqslant 1 + \frac{1}{x} 
	$. 
	
	For simplicity, we denote
	\begin{align}
		B := \frac{2c_1^2A}{R^2} + \frac{c_2 A}{R^2} + \frac{c_1}{R^2}\left(C(N, A)\left(1+R\sqrt{l}\right)+c_0 R\right),
	\end{align}
	which infers that 
	\begin{align}\label{bdB}
		2F^2_{\nabla u}\left(\nabla^{\nabla u}\Psi\right) - \Delta^{\nabla u}\Psi \leqslant B.
	\end{align}
	 {By Cauchy's} inequality and (\ref{cut_F}), it follows that 
	\begin{align}\label{gfH}
		\begin{aligned}
			\Psi g_{\nabla u}\left(\nabla f, \nabla^{\nabla u} H\right) - g_{\nabla u}\left(\nabla f, \nabla^{\nabla u} (\Psi H)\right)&=-g_{\nabla u}(\nabla f, \nabla^{ \nabla u} \Psi) H\\
			& \leqslant F(\nabla f)\left(\frac{c_1}{R}\Psi^{
				\frac{1}{2}}\right) H .
		\end{aligned}
	\end{align}
	 {Since $H$ is defined pointwise, now we can} compute 
	\begin{align}\label{H*Cut}
		\begin{aligned}
			&\Delta^{\nabla u}(\Psi H)=H \Delta^{\nabla u} \Psi+\Psi \Delta^{\nabla u} H+2 g_{\nabla u}\left(\nabla^{\nabla u} \Psi, \nabla^{\nabla u} H\right)  \\
			& =H \Delta^{\nabla u} \Psi+\Psi \Delta^{\nabla u} H-2 H F^2_{\nabla u}\left(\nabla^{\nabla u} \Psi\right) \Psi^{-1} + 2g_{\nabla u}\left(\nabla ^{\nabla u} \Psi, \nabla^{\nabla u} (\Psi H)\right)\Psi^{-1}.\\
		\end{aligned}
	\end{align}
	Using Lemma \ref{lem_f2H}, (\ref{bdB}) and (\ref{gfH}) in (\ref{H*Cut}), a straightforward calculation yields that

	\begin{align}
		\begin{aligned}
			\Delta^{\nabla u}(\Psi H) + 2 g_{\nabla u}\left(\nabla^{\nabla u}\left(f-\log\Psi \right), \mbox{ and } \nabla^{\nabla u} (\Psi H)\right) - \Psi H_t 
			\geqslant  \mathcal{K}, 
		\end{aligned}
	\end{align}
	in sense of distribution on $M \times (0,T)$, where 
	\begin{align}
		\begin{aligned}
			\mathcal{K} :=  t \Psi\left\{\frac{2}{N}\left(-\frac{H}{\beta t}-\left(1-\frac{1}{\beta}\right) F^2(\nabla f)\right)^2+((\beta-1) a-2 K) F^2(\nabla f)\right\}\\
			-B H-\frac{2 c_1}{R} F^2(\nabla f) \Psi^{\frac{1}{2}} H -a \Psi H-\frac{1}{t} \Psi H.
		\end{aligned}
	\end{align}
	Let $(z,s)$ be the maximal point of $\Psi H$ on $B_p(2R) \times[0,T]$ and $\Lambda = \Psi (z)H(z,s)$ be the maximum of $\Psi H$. Following the same argument in section 3, we deduce that $\mathcal{K}(z,s) \leqslant 0$. Here we also set $v := F^2(\nabla f)/H$, and it implies that at $ (z,s)$ 
	\begin{align}
		\begin{aligned}
			& B s(\Psi H)+\frac{2 c_1}{R}(\Psi H)^{\frac{3}{2}} v^{\frac{1}{2}} s+\Psi H+ as\Psi H \\
			& \geqslant \frac{2}{N \beta^2}(\Psi H)^2(1+(\beta-1) v s)^2+((\beta-1) a-2 AK) v s^2(\Psi H) .
		\end{aligned}
	\end{align}
	Note that for any fixed $0 <\delta< 1$,
	\begin{align}
		\begin{aligned}
			\frac{2 c_1}{R}(\Psi H)^{\frac{3}{2}} v^{\frac{1}{2}} s \leqslant&(1-\delta) \frac{2}{\beta^2 N}(\Psi  H)^2(1+(\beta-1) v s)^2 \\
			&+\frac{N \beta^2(\Psi H) v s^2 \cdot \frac{ c_1^2}{R^2}}{2(1-\delta)\left(1+(\beta-1) v s\right)^2},
		\end{aligned}
	\end{align}
	and 
	\begin{align}
		\begin{aligned}
			\frac{N \beta^2 c_1^2 v s^2}{2(1-\delta)(1+(\beta-1) v s)^4 R^2} &\leqslant \frac{N \beta^2 c_1^2 v s^2}{2(1-\delta)(1+4(\beta-1) v s) R^2}\\
			&\leqslant \frac{N \beta^2 c_1^2 s}{8(1-\delta)(\beta-1) R^2}.
		\end{aligned}
	\end{align}
	Therefore we have
	\begin{align}
		\begin{aligned}
			& B s+\frac{N \beta^2 c_1^2}{8(1-\delta)(\beta-1) R^2}s+a s+1 \\
			& \geqslant \frac{2 \delta}{N \beta^2}\Lambda(1+(\beta-1) v s)^2 +((\beta-1) a-2A K) v s^2. \\
		\end{aligned}
	\end{align}
	
	Then by simply modifying how we deal with different cases $a>0$ and $a<0$ in Section 3, one obtains that when $a<0$,
	\begin{align*}
		\Lambda &\leqslant \frac{N \beta^2}{2\delta} \cdot \frac{1+ B s+\frac{N \beta^2 c_1^2}{8(1-\delta)(\beta-1) R^2}s+(-(\beta-1) a +2 AK) s \cdot v s}{\left(1+(\beta-1) v s\right)^2}\\
		&\leqslant\frac{N \beta^2}{2\delta} \left\{ 1+ \frac{B}{2} s+\frac{N \beta^2 c_1^2}{16(1-\delta)(\beta-1) R^2}s - \frac{a}{4}s +  \frac{AK s}{2(\beta-1)}\right\}.
	\end{align*}
	And when $a>0$,
	\begin{align*}
		\Lambda &\leqslant \frac{N \beta^2}{2\delta} \cdot \frac{1+ as + B s+\frac{N \beta^2 c_1^2}{8(1-\delta)(\beta-1) R^2}s +2 AK s \cdot v s}{\left(1+(\beta-1) v s\right)^2}\\
		&\leqslant\frac{N \beta^2}{2\delta} \left\{ 1+\frac{a}{2} s +\frac{B}{2} s+\frac{N \beta^2 c_1^2}{16(1-\delta)(\beta-1) R^2}s  +  \frac{AK s}{2(\beta-1)}\right\}.
	\end{align*}
	
	Finally,
	\begin{align}
		\begin{aligned}
			&\sup_{x\in B_p(R)} \left\{ F^2(\nabla f)(x,T) + \beta af(x,T)  - \beta f_t(x,T) \right\} = \sup_{x\in B_p(R)}\frac{\Psi(x)H(x,T)}{T} \\ 
			&\leqslant\sup_{(x,t)\in B_p(2R)\times [0,T]}\frac{\Psi(x)H(x,t)}{T} = \frac{\Lambda}{T}\\
			&\leqslant \frac{N\beta^2}{2\delta}\left\{\frac{1}{T} + \alpha + \frac{AK}{2(\beta -1)}\right.
			\left.+ \frac{N\beta ^2 c_1^2}{16(1-\delta)(\beta - 1) R^2} + \frac{B}{2} \right\},
		\end{aligned}
	\end{align}
	which concludes the proof.
\end{proof}

\section{Applications}
In this section, we give the two traditional applications of the gradient estimates on Finsler metric measure spaces. We present them in two subsections.
\subsection{Harnack inequalities}
As an application of gradient estimates, We now give Harnack inequalities about the solutions to the equation (\ref{maineq}). Note that by taking  {$a=0$} and $b=0$ the equation (\ref{maineq}) becomes  the Finslerian heat equation and the corresponding results have already been given in  \cite{ohta2014bochner} and \cite{xia2023li}.
\begin{corollary}[Li–Yau-Harnack inequality in compact Finsler  {metric measure space}]\label{cor_hnk_cpt} Under the same conditions in Theorem \ref{mainthm_cpt}, we have
	\begin{align}
		[u(x_1, t_2)]^{1/{e^{a t_1}}} \leqslant[u(x_2, t_2)]^{{1}/{e^{a t_2}}}\exp (\Theta_\beta(x_1,x_2,t_1,t_2)),
	\end{align}
	where 
	\begin{align*}
		\Theta_\beta(x_1,x_2,t_1,t_2) =\int_{t_1}^{t_2} \frac{N\beta}{2e^{a\tau}}\left\{\frac{1}{\tau } + \alpha + \frac{K}{2(\beta -1)} + \frac{d(x_2,x_1)^2}{2N(t_2-t_1)^2}- \frac{2b}{\beta N}\right\} d\tau,
	\end{align*}
	for any $\beta>1$, $0<t_1<t_2<T$, and $x_1$, $x_2\in M$.
\end{corollary}
\begin{proof}
	For the parameter $s\in [t_1,t_2]$ and some tangent vector $v \in T_{x_2}M$, define $\gamma(s):= \exp (x_2,(t_2 - s)v)$ as a reverse curve of the minimal geodesic from $x_2$ to $x_1$. And that implies $\gamma(t_1) = x_1$, $\gamma(t_2) = x_2$ and 
	\begin{align}\label{F_geo}
		F(- \dot{\gamma}) = \frac{d(x_2,x_1)}{t_2-t_1}.
	\end{align}
	Putting 
	\begin{align}
		\begin{aligned}
			& \theta(s)=\frac{N\beta}{2e^{a s}}\left\{\frac{1}{s } + \alpha + \frac{K}{2(\beta -1)} + \frac{d(x_2,x_1)^2}{2N(t_2-t_1)^2}- \frac{2b}{\beta N}\right\},\\
			& \sigma(s)=e^{-a s}f(\gamma(s), s)+\int_{0}^{s}\theta(\tau)d\tau,
		\end{aligned}
	\end{align}
	 {Cauchy's} inequality  {for Finsler metric shows that 
	$$D(f)(\dot{\gamma}) = -g_{\nabla f}(\nabla f, -\dot{\gamma})\geqslant - F(\nabla f)F(-\dot{\gamma})$$
	}
	Then by Theorem \ref{mainthm_cpt},
	and (\ref{F_geo}) we see that
	\begin{align}
		\begin{aligned}
			\sigma^{\prime}(s) &=e^{-as}\left(D f(\dot{\gamma})+f_t-a f \right)+\theta(s)\\
			& \geqslant -e^{-as}F(\nabla f) {F(- \dot{\gamma})} + \frac{F^2(\nabla f)}{\beta e^{as}} + \frac{\beta  d(x_2,x_1)^2}{4e^{as}(t_2-t_1)^2}\\
			& \geqslant - \frac{F(\nabla f)d(x_2,x_1)}{e^{as}(t_2-t_1)} + \frac{F(\nabla f)d(x_2,x_1)}{e^{as}(t_2-t_1)} = 0.
		\end{aligned}
	\end{align}
	Therefore, by the monotonicity of $\sigma$, 
	\begin{align}
		\begin{aligned}
			\exp \left(\frac{f(x_1, t_1)}{e^{a t_1}}\right) \cdot  \exp \left(\int_{0}^{t_1}\theta(\tau)d\tau \right)&=e^{\sigma(t_1)}  \leqslant e^{\sigma(t_2)} \\
			& =\exp \left(\frac{f(x_2, t_2)}{e^{a t_2}}\right) \cdot \exp \left(\int_{0}^{t_2}\theta(\tau)d\tau \right),
		\end{aligned}
	\end{align}
	which demonstrates the assertion.
\end{proof}
\begin{corollary}[Li-Yau-Harnack inequality in non-compact Finsler  {metric measure space}]\label{cor_hnk_cplt} Under the assumption in Theorem \ref{mainthm_cplt}, then we have
	\begin{align}
		[u(x_1, t_2)]^{1/{e^{a t_1}}} \leqslant[u(x_2, t_2)]^{{1}/{e^{a t_2}}}\exp (\Theta_{\beta,\delta}(x_1,x_2,t_1,t_2)),
	\end{align}
	where 
	\begin{align*}
		\Theta_{\beta,\delta}(x_1,x_2,t_1,t_2) =\int_{t_1}^{t_2} \frac{N\beta}{2\delta e^{a\tau}}\left\{\frac{1}{\tau } + \alpha + \frac{AK}{2(\beta -1)} + \frac{\delta d(x_2,x_1)^2}{2N(t_2-t_1)^2}- \frac{2\delta b}{\beta N}\right.\\
		\left.+ \frac{N\beta ^2 c_1^2}{16(1-\delta)(\beta - 1) R^2} + \frac{(2c_1^2 + c_2) A }{2R^2}+\frac{c_3(N,A,K_0)}{2R^2}(1+R+R\sqrt{K}) \right\} d\tau,
	\end{align*}
	for any $0 <\delta< 1$, $\beta>1$, $0<t_1<t_2<T$ and $x_1$, $x_2\in B_p(R)$. 
\end{corollary}
The proof is similar with Corollary \ref{cor_hnk_cpt}. The only different part is to replace the function $\Theta$ related to  {gradient estimate from Theorem \ref{mainthm_cplt}} so we omit the proof here. 
\subsection{ {Boundness and gap of global solutions}}
In this subsection, we turn to the elliptic equation (\ref{e_equa})  { and
apply the aforementioned gradient estimate to prove Theorem \ref{Liou} and Theorem \ref{Lioub0k0}.}

\begin{proof}[Proof  { of Theorem \ref{Liou}}]  Since $u(x,t) \equiv u(x)$, let $R\to + \infty$, $\delta \to 1$ and $t\to \infty$ in (\ref{ge_cplt}), then we have 
	\begin{align}\label{ge_gl}
		F^2(\nabla f)(x) + \beta af(x) + \beta b \leqslant \frac{N\beta^2}{2}\left\{\alpha + \frac{AK}{2(\beta -1)}\right\}.
	\end{align}
	
	The conditions in (\ref{Liou_con}) imply that there exist
	\begin{align}
		\beta_0^\pm :=\frac{4b+2N\alpha -NAK\pm\sqrt{(NAK-4b-2N\alpha)^2-32Nb\alpha}}{4N\alpha} >1.
	\end{align}
	Substituting either of $\beta_0$ 
	 {to (\ref{ge_gl})}, we have
	\begin{align}\label{ge_gl_0}
		F^2(\nabla f)(x) + \beta_0 af(x) \leqslant 0.
	\end{align}
	 {In case of $a<0$, suppose that there exists 
	 an $x_0\in M $ such that $f(x_0)< 0$, then 
	\begin{align*}
		F^2(\nabla f)(x_0) + \beta_0 af(x_0) > 0.
	\end{align*}
	which contradicts (\ref{ge_gl_0}). Hence $\log u(x) =f(x)  \geqslant 0$. The same argument can by applied for the other case.}
\end{proof}
\begin{remark}
	It should be noted that the conditions in (\ref{Liou_con}) are sharp. Failure to fulfill these conditions would result in absence of such a feasible $\beta_0 >1$. Moreover, nonnegative mixed weighted Ricci curvature ($^m\operatorname{Ric}^N\geqslant  0 $) or weighted Ricci curvature ($\operatorname{Ric}^N\geqslant  0 $) inherently meet the second condition in (\ref{Liou_con}).
\end{remark}
%\begin{remark}
%	Since $f = \log u$, one can  derive directly from  {Theorem \ref{liou_thm}} that $u(x)\geqslant 1$ on $M$ when $a<0$. Namely, there is a gap for $u(x)$ when $b>\alpha N/2 = -aN/8$.
%\end{remark}
 {
\begin{proof}[Proof of Theorem \ref{Lioub0k0}]
	Take $b = 0$, and $K=0$ in (\ref{ge_gl}). Noting that $F^2(\nabla f)\geqslant 0$ and $\beta>1$, we have
	\begin{align}
		a f(x)\leqslant \left\{\begin{array}{ll}
			-\dfrac{N\beta a}{8} & \text { when } a<0 ,\\
			\dfrac{N\beta a}{4} & \text { when } a>0.
		\end{array}  
		\right.
	\end{align}
	Then by letting $\beta \to 1$,  we have proved this theorem.
\end{proof}}

%\begin{spacing}{2.0}
%	% 这里是两倍行距
%	\noindent {\textbf {\Large Acknowledgement}}
%\end{spacing}

%The first author is supported partially by the NNSFC (No. 12001099, 12271093).
\section*{Acknowledgments}

The authors are very grateful to the reviewers for their careful review and valuable comments. %The first author is supported partially by the NNSFC (No. 12001099, 12271093).

\section*{Declarations}

\subsection*{Ethical Approval}
Ethical Approval is not applicable to this article as no human or animal studies in this study.

\subsection*{Funding} 
The first author is supported partially by the NNSFC (No. 12001099, 12271093).

\subsection*{Data availability statement}
Data sharing is not applicable to this article as no new data were created or analyzed in this study.
%

%\bibliography{refsab.bib}

\begin{thebibliography}{25}
	\ifx \bisbn   \undefined \def \bisbn  #1{ISBN #1}\fi
	\ifx \binits  \undefined \def \binits#1{#1}\fi
	\ifx \bauthor  \undefined \def \bauthor#1{#1}\fi
	\ifx \batitle  \undefined \def \batitle#1{#1}\fi
	\ifx \bjtitle  \undefined \def \bjtitle#1{#1}\fi
	\ifx \bvolume  \undefined \def \bvolume#1{\textbf{#1}}\fi
	\ifx \byear  \undefined \def \byear#1{#1}\fi
	\ifx \bissue  \undefined \def \bissue#1{#1}\fi
	\ifx \bfpage  \undefined \def \bfpage#1{#1}\fi
	\ifx \blpage  \undefined \def \blpage #1{#1}\fi
	\ifx \burl  \undefined \def \burl#1{\textsf{#1}}\fi
	\ifx \doiurl  \undefined \def \doiurl#1{\url{https://doi.org/#1}}\fi
	\ifx \betal  \undefined \def \betal{\textit{et al.}}\fi
	\ifx \binstitute  \undefined \def \binstitute#1{#1}\fi
	\ifx \binstitutionaled  \undefined \def \binstitutionaled#1{#1}\fi
	\ifx \bctitle  \undefined \def \bctitle#1{#1}\fi
	\ifx \beditor  \undefined \def \beditor#1{#1}\fi
	\ifx \bpublisher  \undefined \def \bpublisher#1{#1}\fi
	\ifx \bbtitle  \undefined \def \bbtitle#1{#1}\fi
	\ifx \bedition  \undefined \def \bedition#1{#1}\fi
	\ifx \bseriesno  \undefined \def \bseriesno#1{#1}\fi
	\ifx \blocation  \undefined \def \blocation#1{#1}\fi
	\ifx \bsertitle  \undefined \def \bsertitle#1{#1}\fi
	\ifx \bsnm \undefined \def \bsnm#1{#1}\fi
	\ifx \bsuffix \undefined \def \bsuffix#1{#1}\fi
	\ifx \bparticle \undefined \def \bparticle#1{#1}\fi
	\ifx \barticle \undefined \def \barticle#1{#1}\fi
	\bibcommenthead
	\ifx \bconfdate \undefined \def \bconfdate #1{#1}\fi
	\ifx \botherref \undefined \def \botherref #1{#1}\fi
	\ifx \url \undefined \def \url#1{\textsf{#1}}\fi
	\ifx \bchapter \undefined \def \bchapter#1{#1}\fi
	\ifx \bbook \undefined \def \bbook#1{#1}\fi
	\ifx \bcomment \undefined \def \bcomment#1{#1}\fi
	\ifx \oauthor \undefined \def \oauthor#1{#1}\fi
	\ifx \citeauthoryear \undefined \def \citeauthoryear#1{#1}\fi
	\ifx \endbibitem  \undefined \def \endbibitem {}\fi
	\ifx \bconflocation  \undefined \def \bconflocation#1{#1}\fi
	\ifx \arxivurl  \undefined \def \arxivurl#1{\textsf{#1}}\fi
	\csname PreBibitemsHook\endcsname
	
	%%% 1
	\bibitem[\protect\citeauthoryear{Gross}{1975}]{gross1975logarithmic}
	\begin{barticle}
		\bauthor{\bsnm{Gross}, \binits{L.}}:
		\batitle{Logarithmic {Sobolev} inequalities}.
		\bjtitle{Am. J. Math.}
		\bvolume{97}(\bissue{4}),
		\bfpage{1061}--\blpage{1083}
		(\byear{1975})
	\end{barticle}
	\endbibitem
	
	%%% 2
	\bibitem[\protect\citeauthoryear{Ohta}{2017}]{ohta2017nonlinear}
	\begin{barticle}
		\bauthor{\bsnm{Ohta}, \binits{S.}}:
		\batitle{Nonlinear geometric analysis on {Finsler} manifolds}.
		\bjtitle{Eur. J. Math.}
		\bvolume{3}(\bissue{4}),
		\bfpage{916}--\blpage{952}
		(\byear{2017})
	\end{barticle}
	\endbibitem
	
	%%% 3
	\bibitem[\protect\citeauthoryear{Yin and Mo}{2021}]{yin2021some}
	\begin{botherref}
		\oauthor{\bsnm{Yin}, \binits{S.}},
		\oauthor{\bsnm{Mo}, \binits{X.}}:
		Some results on complete {Finsler} measure spaces.
		J. Math. Anal. Appl.
		\textbf{497}(1)
		(2021)
	\end{botherref}
	\endbibitem
	
	%%% 4
	\bibitem[\protect\citeauthoryear{Chung and Yau}{1996}]{chung1996logarithmic}
	\begin{barticle}
		\bauthor{\bsnm{Chung}, \binits{F.R.K.}},
		\bauthor{\bsnm{Yau}, \binits{S.-T.}}:
		\batitle{Logarithmic {Harnack} inequalities}.
		\bjtitle{Math. Res. Lett.}
		\bvolume{3}(\bissue{6}),
		\bfpage{793}--\blpage{812}
		(\byear{1996})
	\end{barticle}
	\endbibitem
	
	%%% 5
	\bibitem[\protect\citeauthoryear{Hamilton}{1995}]{hamilton1995formation}
	\begin{barticle}
		\bauthor{\bsnm{Hamilton}, \binits{R.S.}}:
		\batitle{The formations of singularities in the {Ricci Flow}}.
		\bjtitle{J. Differ. Geom.}
		\bvolume{2},
		\bfpage{7}--\blpage{136}
		(\byear{1995})
	\end{barticle}
	\endbibitem
	
	%%% 6
	\bibitem[\protect\citeauthoryear{Zhu}{2022}]{zhu2022class}
	\begin{botherref}
		\oauthor{\bsnm{Zhu}, \binits{H.}}:
		On a class of quasi-{Einstein} {Finsler} metrics.
		J. Geom. Anal.
		\textbf{32}(7)
		(2022)
	\end{botherref}
	\endbibitem
	
	%%% 7
	\bibitem[\protect\citeauthoryear{Ma}{2006}]{ma2006gradient}
	\begin{barticle}
		\bauthor{\bsnm{Ma}, \binits{L.}}:
		\batitle{Gradient estimates for a simple elliptic equation on complete
			non-compact {Riemannian} manifolds}.
		\bjtitle{J. Funct. Anal.}
		\bvolume{241}(\bissue{1}),
		\bfpage{374}--\blpage{382}
		(\byear{2006})
	\end{barticle}
	\endbibitem
	
	%%% 8
	\bibitem[\protect\citeauthoryear{Yang}{2008}]{yang2008gradient}
	\begin{barticle}
		\bauthor{\bsnm{Yang}, \binits{Y.}}:
		\batitle{Gradient estimates for a nonlinear parabolic equation on {Riemannian}
			manifolds}.
		\bjtitle{Proc. Amer. Math. Soc.}
		\bvolume{136}(\bissue{11}),
		\bfpage{4095}--\blpage{4102}
		(\byear{2008})
	\end{barticle}
	\endbibitem
	
	%%% 9
	\bibitem[\protect\citeauthoryear{Huang and Ma}{2010}]{huang2010gradient}
	\begin{barticle}
		\bauthor{\bsnm{Huang}, \binits{G.}},
		\bauthor{\bsnm{Ma}, \binits{B.}}:
		\batitle{Gradient estimates for a nonlinear parabolic equation on {Riemannian}
			manifolds}.
		\bjtitle{Arch. Math.}
		\bvolume{94}(\bissue{3}),
		\bfpage{265}--\blpage{275}
		(\byear{2010})
	\end{barticle}
	\endbibitem
	
	%%% 10
	\bibitem[\protect\citeauthoryear{Garofalo and
		Mondino}{2014}]{Garofalo-Mondino2014}
	\begin{barticle}
		\bauthor{\bsnm{Garofalo}, \binits{N.}},
		\bauthor{\bsnm{Mondino}, \binits{A.}}:
		\batitle{{Li-Yau} and {Harnack} type inequalities in {$RCD^*(K,N)$} metric
			measure spaces}.
		\bjtitle{Nonlinear Anal., Theory Methods Appl., Ser. A, Theory Methods}
		\bvolume{95},
		\bfpage{721}--\blpage{734}
		(\byear{2014})
	\end{barticle}
	\endbibitem
	
	%%% 11
	\bibitem[\protect\citeauthoryear{Fu and Wu}{2022}]{fu2022gradient}
	\begin{barticle}
		\bauthor{\bsnm{Fu}, \binits{X.}},
		\bauthor{\bsnm{Wu}, \binits{J.-Y.}}:
		\batitle{Gradient estimates for a nonlinear parabolic equation with {Dirichlet}
			boundary condition}.
		\bjtitle{Kodai. Math. J.}
		\bvolume{45}(\bissue{1}),
		\bfpage{96}--\blpage{109}
		(\byear{2022})
	\end{barticle}
	\endbibitem
	
	%%% 12
	\bibitem[\protect\citeauthoryear{Dung and Dung}{2019}]{dung2019sharp}
	\begin{barticle}
		\bauthor{\bsnm{Dung}, \binits{H.}},
		\bauthor{\bsnm{Dung}, \binits{N.}}:
		\batitle{Sharp gradient estimates for a heat equation in {Riemannian}
			manifolds}.
		\bjtitle{Proc. Amer. Math. Soc.}
		\bvolume{147}(\bissue{12}),
		\bfpage{5329}--\blpage{5338}
		(\byear{2019})
	\end{barticle}
	\endbibitem
	
	%%% 13
	\bibitem[\protect\citeauthoryear{Ohta and Sturm}{2008}]{Ohta2008HeatFO}
	\begin{barticle}
		\bauthor{\bsnm{Ohta}, \binits{S.}},
		\bauthor{\bsnm{Sturm}, \binits{K.-T.}}:
		\batitle{Heat flow on {Finsler} manifolds}.
		\bjtitle{Commun. Pure Appl. Math.}
		\bvolume{62}(\bissue{10}),
		\bfpage{1386}--\blpage{1433}
		(\byear{2008})
	\end{barticle}
	\endbibitem
	
	%%% 14
	\bibitem[\protect\citeauthoryear{Ohta}{2009}]{ohta2009finsler}
	\begin{barticle}
		\bauthor{\bsnm{Ohta}, \binits{S.}}:
		\batitle{{Finsler} interpolation inequalities}.
		\bjtitle{Calc. Var. Partial Differ. Equ.}
		\bvolume{36}(\bissue{2}),
		\bfpage{211}--\blpage{249}
		(\byear{2009})
	\end{barticle}
	\endbibitem
	
	%%% 15
	\bibitem[\protect\citeauthoryear{Ohta and Sturm}{2014}]{ohta2014bochner}
	\begin{barticle}
		\bauthor{\bsnm{Ohta}, \binits{S.}},
		\bauthor{\bsnm{Sturm}, \binits{K.-T.}}:
		\batitle{Bochner--{W}eitzenb{\"o}ck formula and {Li--Yau} estimates on
			{Finsler} manifolds}.
		\bjtitle{Adv. Math.}
		\bvolume{252},
		\bfpage{429}--\blpage{448}
		(\byear{2014})
	\end{barticle}
	\endbibitem
	
	%%% 16
	\bibitem[\protect\citeauthoryear{Lott and Villani}{2009}]{lott2009ricci}
	\begin{barticle}
		\bauthor{\bsnm{Lott}, \binits{J.}},
		\bauthor{\bsnm{Villani}, \binits{C.}}:
		\batitle{Ricci curvature for metric-measure spaces via optimal transport}.
		\bjtitle{Ann. Math.}
		\bvolume{169}(\bissue{3}),
		\bfpage{903}--\blpage{991}
		(\byear{2009})
	\end{barticle}
	\endbibitem
	
	%%% 17
	\bibitem[\protect\citeauthoryear{Sturm}{2006a}]{sturm2006geometry1}
	\begin{barticle}
		\bauthor{\bsnm{Sturm}, \binits{K.-T.}}:
		\batitle{On the geometry of metric measure spaces. {I}}.
		\bjtitle{Acta Math.}
		\bvolume{196}(\bissue{1}),
		\bfpage{65}--\blpage{131}
		(\byear{2006})
	\end{barticle}
	\endbibitem
	
	%%% 18
	\bibitem[\protect\citeauthoryear{Sturm}{2006b}]{sturm2006geometry2}
	\begin{barticle}
		\bauthor{\bsnm{Sturm}, \binits{K.-T.}}:
		\batitle{On the geometry of metric measure spaces. {II}}.
		\bjtitle{Acta Math.}
		\bvolume{196}(\bissue{1}),
		\bfpage{133}--\blpage{177}
		(\byear{2006})
	\end{barticle}
	\endbibitem
	
	%%% 19
	\bibitem[\protect\citeauthoryear{Xia}{2023}]{xia2023li}
	\begin{barticle}
		\bauthor{\bsnm{Xia}, \binits{Q.}}:
		\batitle{{Li-Yau}’s estimates on {Finsler} manifolds}.
		\bjtitle{J. Geom. Anal.}
		\bvolume{33},
		\bfpage{49}
		(\byear{2023})
	\end{barticle}
	\endbibitem
	
	%%% 20
	\bibitem[\protect\citeauthoryear{Shen}{2023a}]{shenac}
	\begin{botherref}
		\oauthor{\bsnm{Shen}, \binits{B.}}:
		Operators on nonlinear metric measure spaces {I}: A new {L}aplacian comparison
		theorem on {Finsler} manifolds and a traditional approach to gradient
		estimates of {Finslerian} {S}chr\"{o}dinger equation
		(2023)
		{\href{https://arxiv.org/abs/2312.06617}{{arXiv:2312.06617}}}
		{[math.DG]}
	\end{botherref}
	\endbibitem
	
	%%% 21
	\bibitem[\protect\citeauthoryear{Shen}{2023b}]{shenAC2023}
	\begin{botherref}
		\oauthor{\bsnm{Shen}, \binits{B.}}:
		Gradient estimates of the {Finslerian} {Allen-Cahn} equation
		(2023)
		{\href{https://arxiv.org/abs/2401.10894v1}{{arXiv:2401.10894v1}}}
		{[math.DG]}
	\end{botherref}
	\endbibitem
	
	%%% 22
	\bibitem[\protect\citeauthoryear{Cao et~al.}{2013}]{CAO20132312}
	\begin{barticle}
		\bauthor{\bsnm{Cao}, \binits{X.}},
		\bauthor{\bsnm{Ljungberg}, \binits{B.F.}},
		\bauthor{\bsnm{Liu}, \binits{B.}}:
		\batitle{Differential {Harnack} estimates for a nonlinear heat equation}.
		\bjtitle{J. Funct. Anal.}
		\bvolume{265}(\bissue{10}),
		\bfpage{2312}--\blpage{2330}
		(\byear{2013})
	\end{barticle}
	\endbibitem
	
	%%% 23
	\bibitem[\protect\citeauthoryear{Bao et~al.}{2000}]{bao2000introduction}
	\begin{bbook}
		\bauthor{\bsnm{Bao}, \binits{D.}},
		\bauthor{\bsnm{Chern}, \binits{S.-S.}},
		\bauthor{\bsnm{Shen}, \binits{Z.}}:
		\bbtitle{An Introduction to {Riemann-Finsler} Geometry}
		vol. \bseriesno{200}.
		\bpublisher{Springer},
		\blocation{New York}
		(\byear{2000})
	\end{bbook}
	\endbibitem
	
	%%% 24
	\bibitem[\protect\citeauthoryear{Shen}{2001}]{shen2001lectures}
	\begin{bbook}
		\bauthor{\bsnm{Shen}, \binits{Z.}}:
		\bbtitle{Lectures on {Finsler} Geometry}.
		\bpublisher{World Scientific},
		\blocation{Singapore}
		(\byear{2001})
	\end{bbook}
	\endbibitem
	
	%%% 25
	\bibitem[\protect\citeauthoryear{Wu}{2015}]{wu2015comparison}
	\begin{barticle}
		\bauthor{\bsnm{Wu}, \binits{B.}}:
		\batitle{Comparison theorems in {Finsler} geometry with weighted curvature
			bounds and related results}.
		\bjtitle{J. Korean Math. Soc.}
		\bvolume{52}(\bissue{3}),
		\bfpage{603}--\blpage{624}
		(\byear{2015})
	\end{barticle}
	\endbibitem
	
\end{thebibliography}
%\bibliographystyle{abbrv}

\end{document}